\documentclass[twoside]{article}

\usepackage[accepted]{aistats2024}
\usepackage[round]{natbib}
\usepackage[utf8]{inputenc} 
\usepackage[T1]{fontenc}    
\usepackage{hyperref}       
\usepackage{url}            
\usepackage{booktabs}       
\usepackage{amsfonts}       
\usepackage{nicefrac}       
\usepackage{microtype}      
\usepackage{xcolor}         

\hypersetup{
    colorlinks,
    linkcolor={red!25!black},
    citecolor={blue!50!black},
    urlcolor={blue!30!black}
}
\usepackage{transparent}
\usepackage[customcolors, norndcorners]{hf-tikz}

\definecolor{color1}{RGB}{45,45,45}
\definecolor{color2}{RGB}{150,170,160}
\definecolor{color3}{RGB}{220,180,200}
\definecolor{color4}{RGB}{180,50,120}
\definecolor{color5}{RGB}{20,150,200}
\definecolor{color6}{RGB}{200,150,20}
\definecolor{color7}{RGB}{100,100,100}
\definecolor{color8}{RGB}{85,80,140}

\usepackage{amsmath,amsthm,amsfonts,amssymb,amscd}
\usepackage{mathrsfs,bm}
\usepackage{mathtools}

\theoremstyle{definition}
\newtheorem{definition}{Definition} 
\newtheorem{remark}{Remark} 
\usepackage{algorithm, algorithmic}

\usepackage{mathtools} 
\usepackage{booktabs} 
\usepackage{tikz} 
\usepackage{amsmath,amsthm,amsfonts,amssymb,amscd}
\usepackage{mathrsfs,bm}
\usepackage{mathtools}

\newcommand{\ngb}{ u\in\calN^v_{\text{in}}}

\newcommand{\norm}[1]{\|#1\|}

\newcommand{\tlbPsi}{\tilde{\bm{\bPsi}}}
\newcommand{\tlbh}{\tilde{\bh}}
\newcommand{\tlbp}{\tilde{\bp}}
\newcommand{\tlba}{\tilde{\ba}}
\newcommand{\tlbb}{\tilde{\bb}}

\newcommand{\tlbP}{\tilde{\bP}}
\newcommand{\tlbR}{\tilde{\bR}}
\newcommand{\tlbB}{\tilde{\bB}}
\newcommand{\barbR}{\bar{\bR}}

\newcommand{\bLambda}{\bm{\Lambda}}

\newcommand{\bPsi}{\bm{\Psi}}

\newcommand{\bgamma}{\bm{\gamma}}

\newcommand{\bnabla}{\boldsymbol{\nabla}}

\newcommand{\tlbX}{\tilde{\bX}}

\newcommand{\tlbg}{\tilde{\bg}}

\newcommand{\bA}{\mathbf{A}}

\newcommand{\bF}{\mathbf{F}}
\newcommand{\bG}{\mathbf{G}}
\newcommand{\bH}{\mathbf{H}}
\newcommand{\bI}{\mathbf{I}}

\newcommand{\bL}{\mathbf{L}}
\newcommand{\bM}{\mathbf{M}}

\newcommand{\bO}{\mathbf{O}}
\newcommand{\bP}{\mathbf{P}}
\newcommand{\bQ}{\mathbf{Q}}
\newcommand{\bR}{\mathbf{R}}
\newcommand{\bS}{\mathbf{S}}

\newcommand{\bW}{\mathbf{W}}
\newcommand{\bX}{\mathbf{X}}

\newcommand{\bZ}{\mathbf{Z}}

\newcommand{\bU}{\mathbf{U}}
\newcommand{\bC}{\mathbf{C}}
\newcommand{\bB}{\mathbf{B}}

\newcommand{\ba}{\mathbf{a}}
\newcommand{\bb}{\mathbf{b}}
\newcommand{\bc}{\mathbf{c}}

\newcommand{\bg}{\mathbf{g}}
\newcommand{\bh}{\mathbf{h}}

\newcommand{\bn}{\mathbf{n}}
\newcommand{\bp}{\mathbf{p}}

\newcommand{\bx}{\mathbf{x}}

\newcommand{\bz}{\mathbf{z}}

\newcommand{\bbR}{\mathbb{R}}

\newcommand{\calB}{\mathcal{B}}

\newcommand{\calE}{\mathcal{E}}
\newcommand{\calF}{\mathcal{F}}
\newcommand{\calG}{\mathcal{G}}

\newcommand{\calN}{\mathcal{N}}

\newcommand{\calV}{\mathcal{V}}

\newcommand{\calO}{\mathcal{O}}

\newcommand{\barbx}{\bar{\bx}}

\newcommand{\tlbA}{\tilde{\bA}}

\newcommand{\tlbx}{\tilde{\bx}}

\newcommand{\bzero}{\mathbf{0}}
\newcommand{\bone}{\mathbf{1}}

\newcommand{\suml}{\sum\limits}
\newcommand{\minl}{\min\limits}

\newcommand{\bigcapl}{\bigcap\limits}

\newcommand{\tr}{\text{Tr}}

\newcommand{\tth}{\text{th}}

\newcommand{\nwl}{\nonumber\\}

\newtheorem{theorem}{Theorem}
\newtheorem{lemma}{Lemma}

\newcommand{\nf}[1]{\|#1\|_\mathrm{F}} 
\newcommand{\nt}[1]{\|#1\|_2}

\newcommand{\sbjt}{\mbox{{s.t.}}}

\begin{document}

\runningtitle{Asynchronous Decentralized  Optimization with Constraints}

%

\twocolumn[

\aistatstitle{ Asynchronous Decentralized  Optimization with Constraints:  Achievable Speeds of Convergence for Directed Graphs}

\aistatsauthor{ Firooz~Shahriari-Mehr \And Ashkan Panahi }

\aistatsaddress{Chalmers University of Technology \And Chalmers University of Technology} ]

\begin{abstract}

We address a decentralized convex optimization problem, where every agent has its unique local objective function and constraint set. Agents compute at different speeds, and their communication may be delayed and directed. For this setup, we propose an asynchronous double averaging and gradient projection (ASY-DAGP) algorithm. Our algorithm handles difficult scenarios such as message failure, by employing local buffers and utilizing the temporal correlation in the transmitted messages.  We guarantee the convergence speed of our algorithm using performance estimation problems (PEP). In particular, we introduce the concept of the linear quadratic (LQ) PEP. This approach simplifies the analysis of smooth convex optimization problems, going beyond Lyapunov function analyses and avoiding restrictive assumptions such as strong-convexity. Numerical experiments validate the effectiveness of our proposed algorithm.

\end{abstract}

\section{Introduction}
\label{sec:intro}
Consider $M$ computational agents exchanging information over a uni-directional communication network, represented by a directed graph $\calG$. Their goal is to minimize the sum of  local objective functions $f^v(\bx)$ under an intersection of local constraints $S^v$, where $f^v,$ and $S^v$ are only known to node $v$ in $\calG$. This can be written as the following  decentralized constrained optimization problem:
\begin{equation}\tag{P}\label{eq: p1}
    \min_{\bx\in\bbR^{m}} \;\; \frac{1}{M}\suml_{v=1}^{M} f^v(\bx) \qquad \sbjt \quad \bx\in\bigcapl_{v=1}^M S^v.
\end{equation}
In decentralized optimization, each node executes a series of local computations and communicate with its neighboring nodes.
The processing capacity of each node, together with the complexity of the local objective terms, lead  nodes to execute their computations at varying speeds. Communications among agents may also encounter random delays or even message losses,  due to the inherent uncertainties in the underlying communication networks. In this scenario, two different strategies can be considered, known as synchronous and asynchronous.  

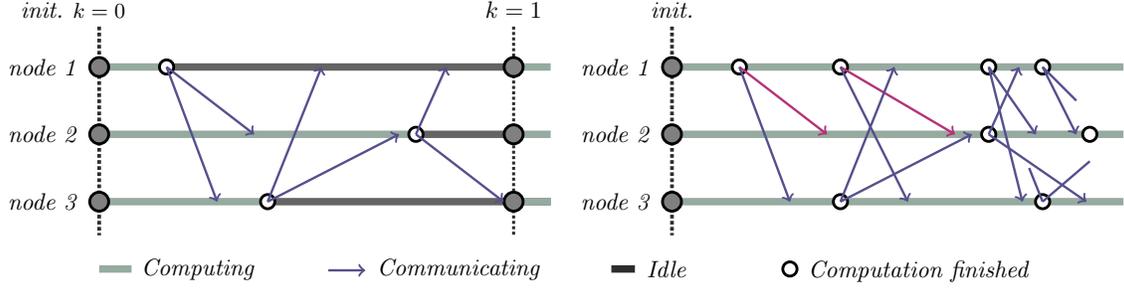
\begin{figure*}[!t]\label{fig: setups}
    \centering
    \resizebox{0.88\linewidth}{!}{\begin{tikzpicture}

\draw[color1, line width=0.55mm, densely dotted, -] (0,0.6) -- (0,-2.5);
\draw[color1, line width=0.45mm, densely dotted] (6.15,0.6) -- (6.15,-2.5);
\draw[color2, line width=1.1mm] (0,0) -- (1,0);
\draw[color7, line width=1.1mm] (1,0) -- (6.15,0);
\draw[color2, line width=1.1mm] (6.15,0) -- (6.7,0);

\draw[color2, line width=1.1mm] (0,-1) -- (4.7,-1);
\draw[color7, line width=1.1mm] (4.7,-1) -- (6.15,-1);
\draw[color2, line width=1.1mm] (6.15,-1) -- (6.7,-1);

\draw[color2, line width=1.1mm] (0,-2) -- (2.5,-2);
\draw[color7, line width=1.1mm] (2.5,-2) -- (6.7,-2);
\draw[color2, line width=1.1mm] (6.15,-2) -- (6.7,-2);

\filldraw[color=black, fill=gray, very thick] (0,0) circle (4pt) node[anchor=west]{};
\filldraw[color=black, fill=gray, very thick] (0,-1) circle (4pt) node[anchor=west]{};
\filldraw[color=black, fill=gray, very thick] (0,-2) circle (4pt) node[anchor=west]{};
\filldraw[color=black, fill=white, very thick] (1,0) circle (3pt) node[anchor=west]{};
\filldraw[color=black, fill=white, very thick] (2.5,-2) circle (3pt) node[anchor=west]{};
\filldraw[color=black, fill=gray, very thick] (6.15,-1) circle (4pt) node[anchor=west]{};
\filldraw[color=black, fill=gray, very thick] (6.15,-2) circle (4pt) node[anchor=west]{};
\filldraw[color=black, fill=gray, very thick] (6.15,0) circle (4pt) node[anchor=west]{};
\filldraw[color=black, fill=white, very thick] (4.7,-1) circle (3pt) node[anchor=west]{};

\draw[->,color8, line width=0.35mm] (1,0) -- (2.3,-1);
\draw[->,color8, line width=0.35mm] (1,0) -- (1.75,-2);
\draw[->,color8, line width=0.35mm] (2.5,-2) -- (3.3,0);
\draw[->,color8, line width=0.35mm] (2.5,-2) -- (4.45,-1);
\draw[->,color8, line width=0.35mm] (4.7,-1) -- (5.15,0);
\draw[->,color8, line width=0.35mm] (4.7,-1) -- (6,-2);

\node [above] at (  0,  0.6) {\footnotesize$k=0$};
\node [above] at (  -0.85,  0.6) {\textit{init.}};
\node [above] at (  6.15,  0.6) {$k=1$};
\node [left]  at (  0,0 )  {\textit{node 1 $\;$} };
\node [left]  at (  0,-1 ) {\textit{node 2 $\;$} };
\node [left]  at (  0,-2 ) {\textit{node 3 $\;$} };




\draw[color1, line width=0.55mm, densely dotted, -] (8.5,0.6) -- (8.5,-2.5);
\draw[color2, line width=1.1mm] (8.5,0) -- (15.2,0);
\draw[color2, line width=1.1mm] (8.5,-1) -- (15.2,-1);
\draw[color2, line width=1.1mm] (8.5,-2) -- (15.2,-2);

\filldraw[color=black, fill=gray, very thick] (8.5,0) circle (4pt) node[anchor=west]{};
\filldraw[color=black, fill=gray, very thick] (8.5,-1) circle (4pt) node[anchor=west]{};
\filldraw[color=black, fill=gray, very thick] (8.5,-2) circle (4pt) node[anchor=west]{};
\filldraw[color=black, fill=white, very thick] (9.5,0) circle (3pt) node[anchor=west]{};
\filldraw[color=black, fill=white, very thick] (11,0) circle (3pt) node[anchor=west]{};
\filldraw[color=black, fill=white, very thick] (11,-2) circle (3pt) node[anchor=west]{};
\filldraw[color=black, fill=white, very thick] (13.2,0) circle (3pt) node[anchor=west]{};
\filldraw[color=black, fill=white, very thick] (14.7,-1) circle (3pt) node[anchor=west]{};
\filldraw[color=black, fill=white, very thick] (13.2,-1) circle (3pt) node[anchor=west]{};
\filldraw[color=black, fill=white, very thick] (14,0) circle (3pt) node[anchor=west]{};
\filldraw[color=black, fill=white, very thick] (14,-2) circle (3pt) node[anchor=west]{};

\draw[->,color4, line width=0.35mm] (9.5,0) -- (10.8,-1);
\draw[->,color8, line width=0.35mm] (9.5,0) -- (10.25,-2);
\draw[->,color4, line width=0.35mm] (11,0) -- (12.7,-1);
\draw[->,color8, line width=0.35mm] (11,0) -- (12,-2);
\draw[->,color8, line width=0.35mm] (11,-2) -- (11.8,0);
\draw[->,color8, line width=0.35mm] (11,-2) -- (12.95,-1);
\draw[-,color8, line width=0.35mm] (14,-2) -- (13.8,-1.5);
\draw[-,color8, line width=0.35mm] (14,-2) -- (14.7,-1.4); 
\draw[-,color8, line width=0.35mm] (14,0) -- (14.5,-0.5);
\draw[->,color8, line width=0.35mm] (14,0) -- (14.5,-1); 
\draw[->,color8, line width=0.35mm] (13.2,-1) -- (13.65,0);
\draw[->,color8, line width=0.35mm] (13.2,-1) -- (14.65,-2);
\draw[->,color8, line width=0.35mm] (13.2,0) -- (13.9,-1);
\draw[->,color8, line width=0.35mm] (13.2,0) -- (13.7,-2);

\node [above] at (  8.5,  0.6) {\textit{init.}};
\node [left]  at (  8.5,0 )  {\textit{node 1 $\;$} };
\node [left]  at (  8.5,-1 ) {\textit{node 2 $\;$} };
\node [left]  at (  8.5,-2 ) {\textit{node 3 $\;$} };

\filldraw[color=black, fill=white, very thick] (10.25,-3) circle (3pt) node[anchor=west]{};
\draw[->,color8, line width=0.35mm] (3.4,-3.02) -- (3.95,-3.02);
\draw[color2, line width=1.1mm] (0,-3) -- (0.48,-3);
\draw[color1, line width=1.1mm] (7.6,-3) -- (7.95,-3);
\node [right] at (  4,  -3) {\textit{Communicating}};
\node [right] at (  10.4,  -3.035) {\textit{Computation finished}};
\node [right] at (  0.5,  -3) {\textit{Computing}};
\node [right] at (  8,  -3) {\textit{Idle}};

\end{tikzpicture}}
    \caption{\textit{Synchronous setup (left):} In each iteration, agents complete their computations and transmit their updated variables to their neighbors. The nodes wait until all messages are delivered to their destination. Then, they simultaneously initiate the next iteration. 
    \textit{Asynchronous setup (right):} Each agent completes its computation and sends the updated variables to its neighbors. Immediately, it starts its next computation using available information. There is a possibility of receiving multiple messages from one node or sometimes no messages. 
    }
    \label{fig: async}
\end{figure*}

Synchronous  optimization methods require all agents to complete their computations and communications before proceeding to the next iteration. 
The speed of these algorithms is limited by the slowest node, called the \emph{straggler}. Additionally, these algorithms require a synchronization clock coordinator, which is a demanding task in large-scale decentralized networks~\citep{hannah2018a2bcd}. Asynchronous decentralized methods 
eliminate the idle time of agents by removing  the synchronization points. In such methods, 
nodes operate uninterruptedly, relying on the current information available from their neighbors.
Figure~\ref{fig: setups} illustrates the distinction between asynchronous and synchronous methods. 

In this paper, we consider the alternative of
asynchronous algorithms 
as they are computationally efficient. We address two major challenges in this setup. 
First, agents must handle delayed information received from their neighbours. 
Second, agents update 
at different rates, which results in difficulties with convergence to the optimal solution of the underlying optimization problem~\citep{assran2020agp}.

\textbf{Algorithmic contribution:}
We introduce, for the first time, an asynchronous algorithm -- the asynchronous double averaging and gradient projection (ASY-DAGP) -- which effectively handles setups with local constraints and directed graphs, ensures guaranteed convergence, and utilizes a fixed step size.
This algorithm is robust to both communication failures and latency, owing to the inclusion of local buffers at each node.
 ASY-DAGP generalizes DAGP, introduced by \citet{firoozTAC}, which is the only synchronous algorithm addressing local constraints over directed graphs, to our knowledge.

\textbf{Theoretical contribution:}
We establish that ASY-DAGP converges  in terms of optimality gap with $O(1/\sqrt{K})$ and its feasibility and consensus gaps vanish with  $O(1/K)$. It is well-known that this rate cannot be improved unless stricter assumptions such as strong convexity are made.
Moreover, we identify the role of delay and asynchrony by a single factor, which we refer to as \emph{delay response}.

In our asynchronous, decentralized, and constrained setup, the standard convergence analysis of optimization algorithms by Lyapunov (potential) functions becomes extremely challenging.
Merely assuming smooth convex functions, we employ a novel proof technique, which draws its inspiration from \citet{firoozTAC}. In particular, we show that \citet{firoozTAC} analysis is a systematic relaxation of the performance estimation problems (PEP) \citep{drori2014performance,taylor2017convex}, which can be used in generic convergence studies. We call this relaxation linear quadratic PEP (LQ-PEP). In asynchronous setups, LQ-PEP has a great advantage that it avoids prerequisite such as bounded delays, being a feature of algorithms e.g. reliant on augmented graphs.

\textbf{Paper organization:} 
In Section \ref{sec: setup and method}, we present our problem setup and detail the ASY-DAGP algorithm. Theoretical results based on LQ-PEP analysis, introduced in Section \ref{sec: LQ-PEP}, are presented in Section \ref{sec: rates}.  Numerical experiments are provided in section~\ref{sec: experiment}.  Finally, we conclude with a literature review in section \ref{sec: litearture}.

\section{Methodology}\label{sec: setup and method}
In this section, we introduce our problem setup and present our algorithm, ASY-DAGP.

\subsection{Problem setup}
We solve the decentralized constrained optimization problem 
 \eqref{eq: p1} in a setup with the following assumptions:
\begin{enumerate}
    \item The local objective functions $f^v$ are convex, differentiable, and $L-$smooth, with $L>0$.
    \item The local constraint sets $S^v$ are closed and convex.
    \item The optimization problem is feasible and there exists a feasible optimal solution $\bx^*$ satisfying the sufficient optimality condition:
    \begin{equation}\label{eq:optimality}
        \bzero \in \sum_{v=1}^M \left( \partial I_{S^v}(\bx^*) + \nabla f^v(\bx^*)\right).
    \end{equation} 
    %
    \item The communication network of the agents is given by a fixed connected directed graph $\calG=(\calV, \calE)$, where the set of vertices $\calV$ represent agents, and $\calE \subseteq \calV\times\calV$ is the set of directed edges. An edge $(v,u)\in\calE$ represents a uni-directional communication link from node $u$ to node $v$. 
    Node  $u$'s incoming and outgoing neighbors are shown by $\calN^u_{\text{in}} = \left\{v \vert (u,v) \in \calE \right\}$ and $\calN^u_{\text{out}} = \left\{ v\vert (v,u)  \in \calE\right\} $, respectively. The graph includes all self-loops; that is,  $u\in\calN^u_{\text{in/out}}$.
    Accordingly, we consider two gossip matrices $\bW$ and $\bQ$, with the same sparsity pattern as the graph adjacency matrix. $\bW$ has zero row sums, and $\bQ$ has zero column sums. They can be built e.g. by the input and output graph Laplacian matrices. 

    \item Communication links are subject to delays. When node $u$ transmits via edge $(v,u)$,  node $v$ receives the message after an arbitrary delay.
    \item 
    Some messages may be lost during communication, i.e. a message transmitted from node $u$ across the communication link $(v,u)$ may not arrive at the intended destination, node $v$.
    \item Agents have varying computational power. Their computations occur at different speeds. As a result, in a given time interval, each node performs a different number of iterations.
    %
\end{enumerate}

DAGP is the only existing synchronous algorithm corresponding to the first four assumptions, but cannot handle the rest. We aim to adapt DAGP for asynchronous setups, addressing the last three  assumptions.

\subsection{ASY-DAGP algorithm}

Similar to DAGP,  ASY-DAGP attempts to solve \eqref{eq:optimality} in an iterative and distributed fashion, by splitting, i.e. re-writing it as:
\begin{equation}\label{eq:split}
          \sum_{v=1}^M\bg^v=\bzero,\quad \bg^v\in\partial I_{S^v}(\bx^*) + \nabla f^v(\bx^*).
    \end{equation} 
the relation $\sum_{v=1}^M\bg^v=\bzero$ is called \emph{null condition}. 
In DAGP, null condition is satisfied by the \emph{distributed null projection} mechanism that uses another set of auxiliary variables $\bh^v$, always satisfying the null condition by design. The variable $\bg^v$ tracks $\bh^v$, that is upon convergence $\bg^v$ becomes equal to $\bh^v$, hence satisfying the null condition.  For comparison, DAGP is detailed in Algorithm \ref{alg: dagp}.

In ASY-DAGP, we attempt to perform a similar procedure as in DAGP. However, the asynchrony affects many steps, for which we introduce suitable modifications. ASY-DAGP is presented in Algorithm \ref{alg: asy-dagp}.
In this algorithm, agents update their local variables in parallel with each other. When an agent executes an iteration and sends the updated variables to its neighbors, it instantly starts its next iteration. However, the duration of an iteration can be variable.
We introduce local buffers, \(\calB^{vu}\), to store all messages node \(v\) receives from \(u\). After node \(v\) processes its messages, the buffer is cleared.

Although each node broadcasts the same message to all its out neighbours, they receive and process that message at different iterations due to asynchrony. This leads to a variable number of messages in a buffer. When there are one or more messages in a buffer $\calB^{vu}$, two averages $\ba^{vu}$ and $\bb^{vu}$ of them are calculated by \eqref{eq: estimates_a} and \eqref{eq: estimates_b}. Otherwise, the last calculated averages are reused. Comparing \eqref{eq: z_update} and \eqref{eq: asy-z_update}, as well as \eqref{eq: h_update} and \eqref{eq: asy-h_update},  we observe that $(\ba^{vu},\bb^{vu})$ substitutes the message pair $(\bx^v,\bh^v-\bg^v)$ in DAGP. Our idea is that due to temporal correlation in the transmitted messages, the former tracks the latter and in effect ASY-DAGP functions similar to DAGP.

There is however a difficulty with the above strategy. Since, each node calculates a distinct average of their buffer, in general $\bb^{vu}\neq\bb^{wu}$, for some $ v\neq w$. A consequence of this phenomenon is that the variables $\bh^v$ will not always satisfy the null condition, and 
the distributed null projection strategy fails. To address this, we introduce new vectors $\bp^v$ and modify the update rule of \eqref{eq: h_update}. In appendix \ref{sec: fixed-point}, we show that convergence to a consensus solution necessitates $\bg^v=\bh^v=\frac{1}{\gamma-1}\sum_{u\in\calN^v_{\mathrm{in}}}q_{vu}\bp^u$. Since $\bQ$ has zero column sums, null condition is implied. As seen, ASY-DAGP requires an additional variable $\bp^v$, but an advantage observed by the analysis of the fixed point in appendix \ref{sec: fixed-point} is that the condition $\ker(\bQ) = \ker(\bW^T)$ in DAGP is no longer necessary. 
\begin{algorithm}[t!]
\caption{DAGP {\footnotesize\citep{firoozTAC}.}}
\label{alg: dagp}
	\begin{algorithmic}[1]
		\renewcommand{\algorithmicrequire}{\textbf{Input:}}
		\REQUIRE step size $\mu$, scaling parameters $\alpha$ and $\rho$, and gossip matrices $\bW$ and $\bQ$.
		\STATE Initialize $k=0$. Initialize $\bx^v_0$ randomly, and
             $\bg^v_0$ and $\bh^v_0$ with zero vectors, \hspace{0.02cm} $\forall v\in\calV$.
		\REPEAT 
		\STATE Update $\bz^v, \bx^v, \bg^v$, and $\bh^v$ variables, $\;\forall v\in\calV$:
              \begin{alignat}{3}
                &\bz^v_{k+1} && = \; && \bx^v_k - \suml_{\ngb} w_{v u}\bx^u_k - \mu \left( \nabla f^v(\bx^v_k) -  \bg^v_k   \right) \label{eq: z_update}\\
                &\bx_{k+1}^v && = \; &&\mathrm{P}_{S^v} \left( \bz^v_{k+1} \right) \label{eq: x_update} \\
                &\bg_{k+1}^v && = \; &&\bg_k^v + \rho (\nabla f^v(\bx^v_k)  - \bg^v_k ) \nonumber\\ 
                & && &&  \;  + \frac{\rho}{\mu}\left(\bz^v_{k+1} - \bx^v_{k+1} \right)    + \alpha \left( \bh_k^v - \bg_k^v \right) \label{eq: g_update} \\
                &\bh_{k+1}^v && = \; &&\bh_k^v -\suml_{\ngb} q_{v u}(\bh_k^u - \bg_k^u) \label{eq: h_update}
            \end{alignat}
		\STATE Send the updated tuple $(\bx^v_{k+1}, \; \bh^v_{k+1}- \bg^v_{k+1})$ to all out-neighbors $u\in\calN^v_{\mathrm{out}}$,  $\;\forall v\in\calV$.
            \STATE Update iteration index: $k = k+1$.
		\UNTIL{Convergence}
	\end{algorithmic}
\end{algorithm}

\begin{algorithm}[t!]
\caption{ASY-DAGP.}
\label{alg: asy-dagp}
	\begin{algorithmic}[1]
		\renewcommand{\algorithmicrequire}{\textbf{Input:}}
		\REQUIRE step size $\mu$, scaling parameters $\alpha$, $\rho$, $\eta$,  gossip matrices $\bW$ and $\bQ$, and forgetting factor $\gamma$.
		\STATE Initialize $\bx^v$ randomly. Initialize $\bg^v$, $\bh^v$, and $\bp^v$ with zero vectors, $\;\forall v\in\calV.$
            \STATE Initialize $\bb^{vu} = \bzero$, $\forall v,u\in\calV$,  $\ba^{vv} = \bx^v$, $\forall v\in\calV$, and $\ba^{vu} = \bzero$, $\forall v,u\in\calV, u \neq v.$
            \STATE Agents continuously receive information from their neighbors and store it in their buffers $\calB^{vu}$.
            \STATE Agents update their variables in parallel as: 
		\REPEAT 
             \STATE Update $\bz^v, \bp^v,$ and $\bh^v$ variables, 
            \begin{alignat}{3}
                &\bz^v && = \; && \bx^v - \suml_{\ngb} w_{v u}\ba^{vu} - \mu \left( \nabla f^v(\bx^v) -  \bg^v   \right) \label{eq: asy-z_update}\\
                &\bp^v && = \; && \bp^v - \eta\suml_{\ngb} q_{vu}\bb^{vu} + \eta(\gamma-1)\bg^v \label{eq: asy-p_update}
                \\
                &\bh^v && = \; && \gamma\bh^v -\suml_{\ngb} q_{v u}\bb^{vu} \label{eq: asy-h_update}
            \end{alignat}
                 \STATE Update $\bx^v$, $\bg^v$ variables similar to \eqref{eq: x_update}-\eqref{eq: g_update}.
    		\STATE Send the updated tuple $(\bx^v, \; \bp^v)$ to all out-neighbors in $\calN^v_{\mathrm{out}}$.
            \STATE Update the estimates of $\bx^u$ and $\bp^u$  at node $v$ for every $u\in\calN^v_{\text{in}}$.
            If the buffer $\calB^{vu}$ is empty, reuse the existing value in memory; otherwise, compute the new estimates as:
            \begin{equation} \label{eq: estimates_a}
                \ba^{vu} = \frac{1}{|\calB^{vu}|} \suml_{\bx^u\in \calB^{vu}} \bx^u  
            \end{equation}
            \begin{equation} \label{eq: estimates_b}
                \bb^{vu} = 
                     \frac{1}{|\calB^{vu}|} \suml_{\bp^u\in \calB^{vu}} \bp^u
            \end{equation}
            \STATE Clear the buffers $\calB^{vu}$, for every $ u\in\calN^v_{\text{in}}$.
		\UNTIL{Convergence}
	\end{algorithmic}
\end{algorithm}

\section{Convergence guarantees}\label{sec: rates}

For the sake of analysis, we index the iterations of each node, individually. Hence, iteration $k$ at different nodes happen at different real times. However, it is harmless to pretend that each iteration occurs simultaneously at every node, but the messages experience a corrected delay that may be negative (non-causal). This assumption does not affect the generality of our approach. 
Based on this convention, we introduce the local index sets $T^{vu}_k$ containing all iteration numbers $k^\prime$, where a message from $u$ is sent  at $k^\prime$ and is processed by $v$ at time $k+1$. In other words, $T^{vu}_k = \{k^\prime | k^\prime+\tau^{vu}_{k^\prime}=k+1\}$.

 Further, we need to introduce few definitions.  First,  we formalize a concept, originally employed by \citet{firoozTAC}, which incorporates frequency analysis tools into the study of optimization algorithms.

\begin{definition}[Proper matrices]
    A tuple $(\bS,\bar\bR,\bP)$ is  \emph{proper} if, for any given complex value $z$ and a positive real value $\beta$, all the statements from Definition 5 in \citep{firoozTAC} hold for the following \emph{forward-backward transfer matrix}:
    \begin{equation}\label{eq: F_define}
    \bF_{\beta}(z) = \left[ 
    \begin{array}{cc}
         \bS & \bI-z^{-1}\barbR^T \\
         \bI-z\barbR & \frac{1}{\beta}\bP\bP^T
    \end{array}
    \right]
    \end{equation}
In particular, we denote by $z_i(\beta)$ the simple roots of $\det(\bF_\beta(z))$.
\end{definition}

Next, we introduced a parameter quantifying the impact of delays and asynchrony.

\begin{definition}[Delay response]
    For a nonempty index set $T_k^{vu}$, define
    \begin{equation}\label{eq: 133}
        \tau^{vu}_{i, k}(\beta)=\frac 1{|T_k^{vu}|}\suml_{k^\prime\in T_k^{vu}} z_i^{|k+1-k^\prime|}(\beta)-1
    \end{equation}
    as the $i^\tth$ delay power spectrum of link $(v,u)$ at time $k$. Accordingly, we define the  delay response of link $(v,u)$ at time $m$ as
    \begin{equation}\label{eq: 144}
        \kappa^{vu}_{m}(\beta)=\suml_{k\mid T_k^{vu}\neq\emptyset, i}|z_i(\beta)|^{|k-m|}\tau^{vu}_{i, k}(\beta),
    \end{equation}
    Then, the delay response $\kappa$ is defined as the supremum of $\kappa^{vu}_{m}(\beta)$ over every time $m$ and link $(v,u)$ and $\beta$.
\end{definition}
Our main result is based on three matrices $\bS,\bar\bR,\bP$ which only depend on the parameters of the algorithm and the network topology. Due to space limitation, they are given in Appendix \ref{sec: convergence}.

\begin{theorem}[Rates of convergence] \label{thm: theorem1}
    Suppose that $(\bS,\bar\bR,\bP)$ is proper. 
    There exist a  constant $C$ independent of the delays such that for $\kappa<\frac 1C$ and a sufficiently large $K$, ASY-DAGP converges as follows: 
    \begin{equation}  
    \nt{\bar\bx^v_K - \bar\bx_K}^2 \leq \frac{CC_0}{ MK},
    \end{equation}
    \begin{equation}
    \left|\suml_v \left( f^v(\barbx^v_K)- f^v(\bx^*)\right)\right|\leq 
    \frac{CC_0}{\mu K}+\sqrt{\frac{C_2CC_0}{ M K}},
    \end{equation}
    where $K$ indicates that all nodes have performed $K$ iterations,  $\bar\bx^v_K = \frac{1}{K}\sum_{k=0}^{K-1}\bx^v_k$,  $\bar\bx_K = \frac{1}{M}\sum_{v=1}^M \bar\bx^v_K$. $C_0$ is a constant dependent on the initial point, and $C_2= \sqrt{\sum_v\|\bn^v+\nabla f^v(\bx^*)\|^2}$,  measuring the heterogeneity of local functions.  
\end{theorem}

%

\begin{remark}
     Note that by $K$ we mean that all nodes have performed $K$ iterations.  However, asynchrony could cause nodes to complete the $K\tth$ iteration at various times. As a result, the convergence in real time hinges on the last node completing its $K\tth$ iteration. 
\end{remark}

\begin{remark}
The assumption that the tuple of matrices is proper is discussed in Remark 2 of \citep{firoozTAC}, where a similar transfer matrix is studied. This study indicates that being proper matrices is not a restrictive assumption, but dictates the range of parameters yielding convergence. 
\end{remark}


In the following section, we provide a sketch of the proof. The detailed proof, the description of the matrices, along with their structures, are provided in Appendix~\ref{sec: convergence}.

\section{Proof outline of Theorem \ref{thm: theorem1}}\label{sec: LQ-PEP}
A standard approach for studying the convergence speed of optimization algorithms is  the Lyapunov-based analysis~\citep{polyak1987introduction}, discussed in many classical papers, e.g.~\citep{nesterov2012efficiency, defazio2014saga,schmidt2017minimizing}. This method relies on finding a suitable Lyapunov (potential) function, which becomes difficult as the algorithm becomes complex. This is especially the case for distributed and asynchronous optimization methods. The  DAGP analysis in \cite{firoozTAC} is based on an alternative approach, called aggregate lower bounding (ALB). 
In the following, we show the ALB is a relaxation of the celebrated performance estimation problem (PEP), introduced by \citet{drori2014performance}, which is a systematic ways to study optimization algorithms. We refer to the relaxation underlying ALB as LQ-PEP, based on which we provide convergence rates for ASY-DAGP. 

LQ-PEP offers several advantages. It incorporates frequency analysis tools, particularly Fourier analysis, into the study of optimization algorithms. Additionally, LQ-PEP provides a tighter convergence bound than Lyapunov analysis when employing the same inequalities in the analysis process.

\subsection{Performance Estimation Problem} 
Let us first take a brief look at PEP. It evaluates the worst-case performance of optimization algorithms after a fixed number $K$ iterations.  An algorithm consists of a state vector $\bPsi_k$ evolving over different iterations $k=0,1,\ldots,K-1$, based on the underlying objective function $f$, and concludes by returning a terminal state $\bar\bPsi_{K-1}$.  %
Note that we make a distinction between a state $\bPsi$ and the solution $\bX$. For our result in Theorem \ref{thm: theorem1}, a state consists of all variables of ASY-DAGP at every node and the terminal state is the average state $\bar\bPsi_{K-1}=\frac 1K\sum_{k=1}^K\bPsi_K$.
The performance is quantified by a performance measure $\Phi(\bar\bPsi_{K-1}, \bPsi_*)$, which is a positive-definite function. For us, this is a combination of terms corresponding to optimality gap, feasibility gap and consensus.
PEP is defined as a search problem for the worst performance over all objective functions $f$ from a given family $\calF$: 
\begin{align}
    \max_{f, \bPsi_*, \{\bPsi_k\}_{k=0}^{K-1}} \quad & \Phi(\bar\bPsi_{K-1}, \bPsi_*) \label{eq: PEP} \\ 
    \sbjt \quad & f \ \text{be a function in family} \ \calF,  \nonumber \\ 
    & \bPsi_* \ \text{is the minimizer of} \ f, \label{cons_pep2}  \nonumber \\
    & \| \bPsi_0 - \bPsi_* \| \leq R,   \nonumber\\
    & \text{algorithm generates} \ \bPsi_1,\ldots,\bPsi_{K-1},   \nonumber
\end{align}
where the dependence on $f$ often makes the problem infinite-dimensional. 
\citet{taylor2017smooth} and \citet{taylor2017convex} make an interesting observation that the PEP in \eqref{eq: PEP} can be reduced to a finite-dimensional problem (finite PEP). 
The finite PEP is a complicated non-convex optimization.
Therefore, its convex relaxed versions are solved by \citet{drori2014performance} and \citet{ taylor2017smooth}. Solving the relaxed finite PEPs
naturally leads to  upper bounds on the speed of convergence, without any need for designing Lyaponov functions.

\subsection{Linear Quadratic (LQ-)PEP}

In \cite{firoozTAC}, the ALB approach is introduced, which amounts to a novel finite-dimensional relaxation of \eqref{eq: PEP} that we call LQ-PEP. This relaxation is obtained in two steps of a) finding a quadratic upper bound on the performance metric $\Phi(\bar\bPsi_{K-1},\bPsi_*)$ and b) a linear relaxation on the constraints in \eqref{eq: PEP}. In Appendix \ref{sec: convergence}, we show that the upper bound on the performance metric is given by 
\begin{equation}\label{eq:quadratic}
    \Phi(\bar\bPsi_{K-1},\bPsi_*)\leq -\frac 1K \sum_{k=0}^{K-1} \langle \tilde\bPsi_k,\bS\tilde\bPsi_k\rangle,
\end{equation} where $\bS$ is a fixed matrix depending on the topology of the network.
We also show that the following relation holds among the states generated by the algorithm:
\begin{equation}\label{eq:linear}
    \tlbPsi_{k+1} - \bar\bR\tlbPsi_k -  \bP\tilde\bU_{k} = \suml_{l=0}^{K-2} \tilde\bR_{k,k+1-l}\tlbPsi_l,
\end{equation}
where for simplicity we introduce $\tlbPsi_k \coloneqq \bPsi_k- \bPsi_*$ for a suitable optimal state $\bPsi_*$, and $\tilde\bU_k$ is a sequence calculated by the states $\tlbPsi_{k'}$. Here, $\bP$ is a fixed matrix, $\barbR$ only depends on the parameters of the algorithm and $\tlbR_{k,\tau}$ only depend on the delay pattern. Without any delay, we have $\tlbR_{k,\tau}=\bzero$.

%

Then, LQ-PEP is given by the following result of the above statements:
\newtheorem{prop}{Proposition}
\begin{prop}\label{prop1}
Take matrices $\bS,\barbR,\bP$ and $\{\tilde\bR_{k,\tau}\}$  as in \eqref{eq:linear}, and a performance measure $\Phi$ satisfying \eqref{eq:quadratic} and encompassing the optimality gap, the feasibility gap and the consensus error.
    Then,  the relation $\Phi(\bar\bPsi_{K-1},\bPsi_*)\leq\frac{C\|\tilde\bPsi_0\|^2}{K}$ holds if the optimal value of the following linear quadratic optimization is zero. 
    \begin{gather}
        \minl_{\substack{\{\tilde\bPsi_k\}_{k=0}^{K-1}, \\ {\{\tilde\bU_{k}\}_{k=0}^{K-2}}}}
        \quad \frac{1}{2}\suml_{k=0}^{K-1}\langle\tilde\bPsi_k,\bS\tilde\bPsi_k\rangle+\frac{C}{2}\|\tlbPsi_0\|^2_{\mathrm{F}} \label{eq: ALB_opt}\\ 
        \text{subject to} \nonumber \\
        \tlbPsi_{k+1} - \bar\bR\tlbPsi_k -  \bP\tilde\bU_{k} = \suml_{l=0}^{K-2} \tilde\bR_{k,k+1-l}\tlbPsi_l.  \nonumber \\
        k=0,1,\ldots, K-2  \nonumber
    \end{gather}
\end{prop}
The explicit construction of $\bS,\bar\bR,\tilde\bR_{k,l}$ and $\Phi$ is given in Appendix \ref{sec: convergence}. Note that $\bS$  is an indefinite matrix and the condition in \eqref{eq: ALB_opt} may only be achieved by considering the constraints. The constraints in \eqref{eq: ALB_opt} represent a linear non-causal dynamical system, with $\tilde\bU_k 
$ as its input. If there are no delays, 
the dynamic will become time invariant. Moreover, in this system, the existence of negative delays is considered in the right hand side of the dynamics.

The proof of proposition \eqref{prop1} is straightforward. Suppose that the optimal value is zero. Note that any sequence of states generated by the algorithm is feasible and hence satisfies $\frac{1}{2}\sum_{k=0}^{K-1}\langle\tilde\bPsi_k,\bS\tilde\bPsi_k\rangle+\frac{C}{2}\|\tlbPsi_0\|^2_{\mathrm{F}}\geq 0$. By \eqref{eq:quadratic}, we see that the result of proposition \eqref{prop1} holds.

\subsection{Establishing Theorem \ref{thm: theorem1}}
To obtain Theorem \ref{thm: theorem1}, we solve \eqref{eq: ALB_opt} and show that the conditions of LQ-PEP in Proposition \ref{prop1} hold. Here, we follow \cite{firoozTAC} and use the Lagrangian method of   multipliers to show that this condition is equivalent to the following system of recurrences having no non-zero solution:
\begin{align}
    \tlbPsi_{k+1} - \bar\bR\tlbPsi_k - \frac{1}{\beta} \bP\bP^T\bLambda_k &= \suml_{l=0}^{K-2} \tilde\bR_{k,k+1-l}\tlbPsi_l \nonumber \\
    \bLambda_{k-1} - \bar\bR^T\bLambda_k + \bS\tlbPsi_k &= \suml_{l=0}^{K-2} \tilde\bR^T_{l,l-k+1}\bLambda_l, \nonumber
\end{align}
where $\bLambda_k$s are the Lagrangian multipliers.  The proof of this fact and the rest of the proof of Theorem \ref{thm: theorem1} is given in Appendix \ref{sec: lqpep}.

\section{Experimental results}\label{sec: experiment}

\begin{figure*}[!t]
\centering
\begin{tabular}{cccc}
\includegraphics[width=0.41\textwidth]{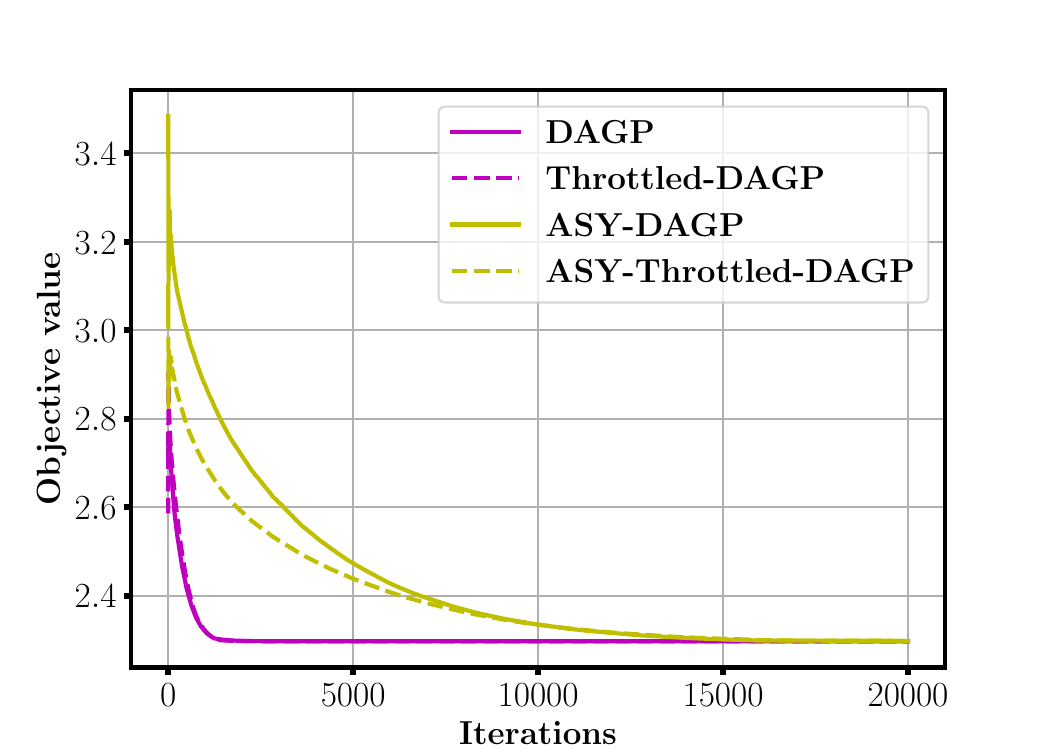} &
\unskip\ \vrule\
\includegraphics[width=0.41\textwidth]{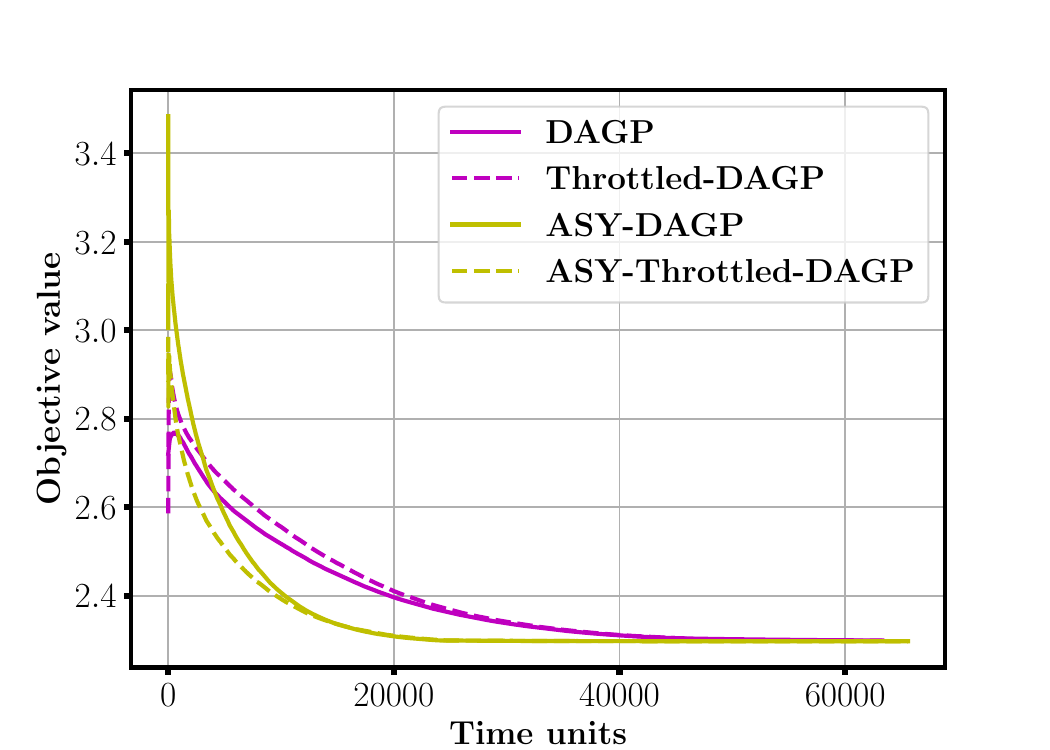} \\ [2pt]
\textit{(a) Convergence over iterations.} & \textit{(b) Convergence over time.}
\end{tabular}
\begin{tabular}{cccc}
\includegraphics[width=0.41\textwidth]{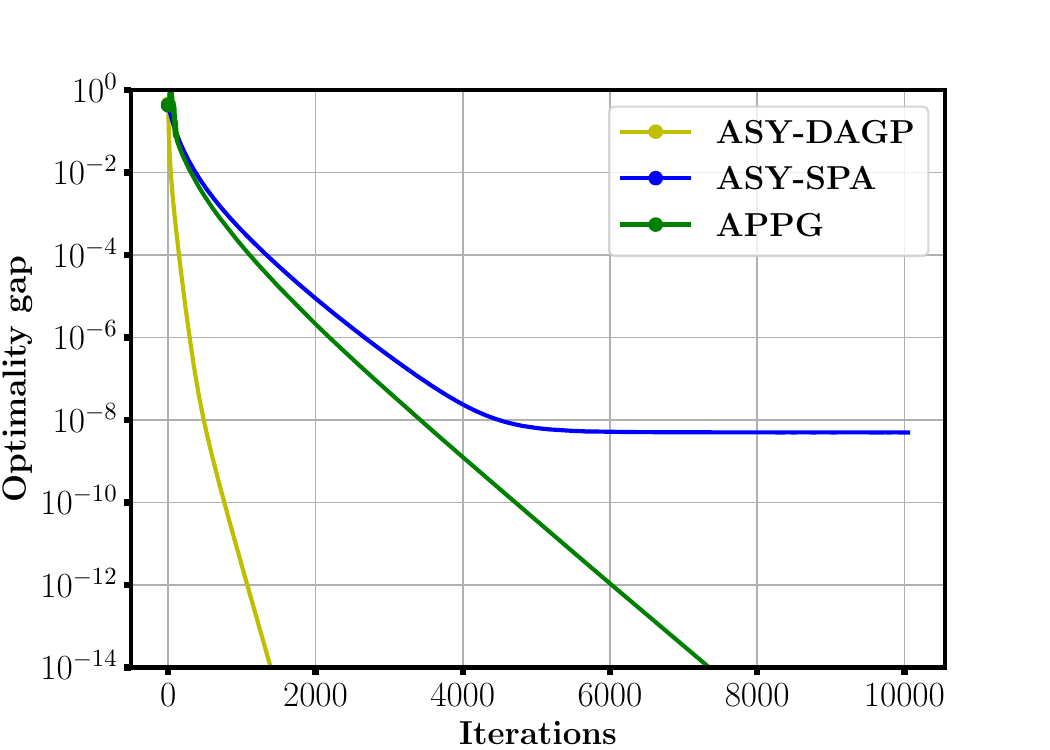} &
\unskip\ \vrule\
\includegraphics[width=0.41\textwidth]{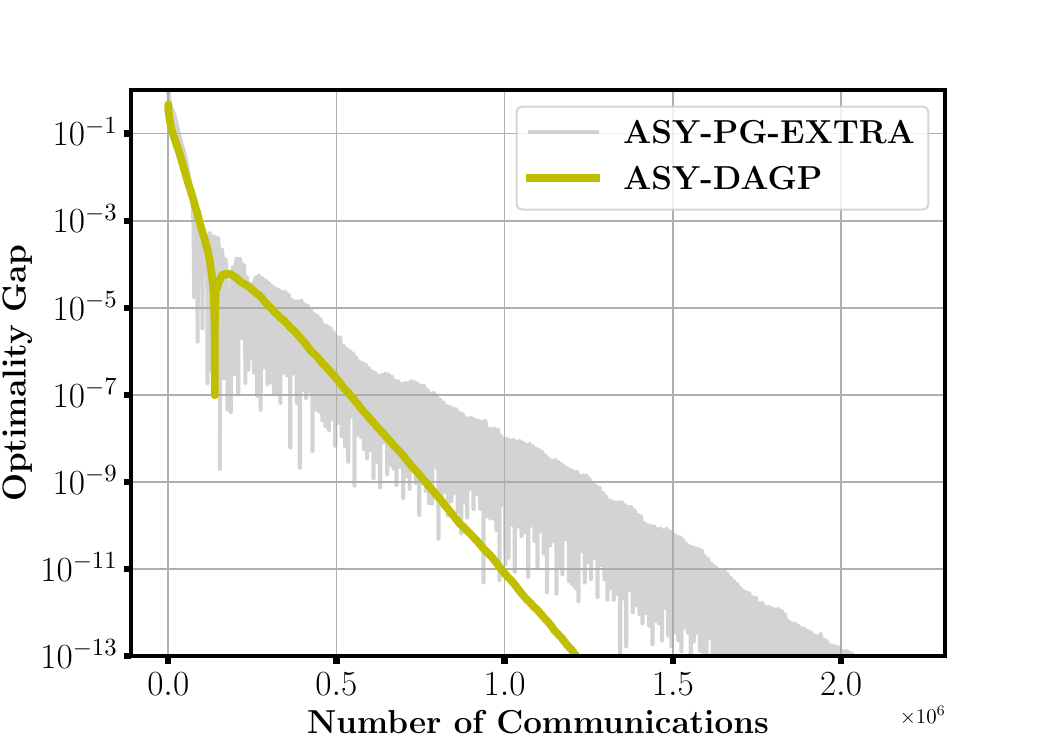} \\ [2pt]
\textit{(c) Unconstrained setup.} &\textit{ (d) Setup with undirected graph. }

\end{tabular}
\caption{ 
\textit{(a,b)} Solving a constrained problem, and drawing a comparison to DAGP and its throttled version.
\textit{(c)} Solving unconstrained logistic regression problem, and comparing to APPG and ASY-SPA. 
\textit{(d)} Solving a constrained problem over undirected graphs, and comparing to ASY-PG-EXTRA. 
}
\label{fig: all}
\end{figure*}

In our experiments, we utilized  directed and random networks, with edge probability $p_e$. We assumed that the  delays for each communication link follow an exponential distribution with mean $\tau_{\text{comm}}^{vu}$. In the first experiment, the computation time for each node followed a uniform distribution within $[1, \bar\tau_{\text{comp}}^{v}]$. In the second one, they were exponentially distribution with  mean $\tau_{\text{comp}}^v$. We used scaled Graph Laplacian matrices as gossip matrices, and the design parameters were manually tuned for all algorithms. Each experiment was conducted multiple times, but we present only a single instance due to the negligible variations across separate runs. We report the objective values computed at the average solution across all nodes.
In our simulations, the x-axis shows the global iteration counter, which increases when at least one agent computes. 

In every experiment, we simulated a distributed setup on a single machine by sequentially performing the actions of different nodes. This approach gives better control over parameters and facilitates the implementation of asynchrony, delays, and message failures. As our experiments have an illustrative nature, we postpone more realistic implementations to a future study.


To the best of our knowledge, there is no other algorithm like ASY-DAGP that handles  a constrained optimization with asynchrony over directed graphs. Hence, we compare ASY-DAGP to other existing algorithms in  simplified  setups, which still demonstrate ASY-DAGP's advantages. Specifically:
\begin{itemize}
\item In the first experiment, we use a directed and constrained setup, but it is synchronous.
\item In the second experiment, our setup is asynchronous and directed, but without constraints.
\item The third experiment features an asynchronous setup that is constrained, but undirected.
\end{itemize}

\textbf{First experiment:}
In this experiment, we compare ASY-DAGP with DAGP. We consider the following local objective functions and constraints
\begin{gather}
    f^v(\bx)  = \log\big(\cosh( \ba_v^T\bx - b_v)\big), \label{eq: synthetic_objective} \\
S^v = \{ \bx \; | \;  \bc_v^T\bx - d_v \leq 0 \},
\label{eq: synthetic_constraint}
\end{gather}
where constants and coefficients are generated from normal distributions. We initialize both algorithms randomly and assume  $M=10, \bx\in\bbR^{5}, p_e=0.8, \tau^{vu}_{\text{comp}}=10, \bar\tau^v_{\text{comp}} = 5v,$ where $v$ indicates the node number. We set the design parameters to 
$\rho       = 0.01,         
\alpha     = 0.1,
\gamma     = 0.5,
\eta       = 1.0.$
Figure \ref{fig: all}a and \ref{fig: all}b present the resulting objective values. We observe both algorithms converge to the same point, demonstrating that ASY-DAGP is capable of reaching the optimal solution of \eqref{eq: p1}. In terms of convergence over iterations, ASY-DAGP requires more iterations. However, it achieves significantly faster wall-clock convergence. It is worth noting that an increase in heterogeneity of the local computation time (called throttling) - for instance, by deliberately slowing down several nodes by a factor of two - may result in a slower convergence for the synchronous setup. For ASY-DAGP,  this situation is not harmful.

\textbf{Second experiment:}
We compare ASY-DAGP to the APPG algorithm \citep{zhang2019fully} in an asynchronous setup. We consider an unconstrained logistic regression problem with $\ell_2$ regularization, similar to the one presented in \citep{firoozTAC}. This problem is applied to two digits of the MNIST dataset~\citep{lecun2010mnist}. A random graph with $p_e=0.6$ and $M=20$ nodes is considered. We use $N_s=1000$ samples from the dataset and a regularization factor of $1/N_s$.
The optimal solution is computed by running the centralized gradient descent for a substantial number of iterations. The optimality gap, i.e., the error in the objective function, is reported in Figure~\ref{fig: all}c.
In this experiment, the design parameters are set to 
$
\tau^{vu}_{\text{comm}} = 50,
\bar\tau^v_{\text{comp}} = 5v,
\gamma     = 0.5,
\eta       = 1.0,
\rho       = 0.1,
\alpha     = 0.7,
\mu        = 1.0.
$
For both APPG and ASY-DAGP, we selected similar step-sizes slightly smaller than $2/L$, provided they do not diverge. Although ASY-SPA typically requires a diminishing step-size, we chose a small fixed one. With this choice, we may only arrive at a neighborhood of the optimal solution.

\textbf{Third experiment:}
In this experiment, we use the same setup as in the first experiment, with the key difference of utilizing an undirected graph. We compare ASY-DAGP with the ASY-PG-EXTRA algorithm in \citep{sayed2017decentralized}. For both algorithms, we employ a near optimal step-size, and for the ASY-PG-EXTRA, we set the relaxing parameter to $0.8$ for all agents. The results are shown in Figure \ref{fig: all}d. 
We plot the optimality gap with respect to the number of  completed communications. ASY-DAGP converges to the optimal solution with only half the communications needed by ASY-PG-EXTRA to achieve the same result.

\begin{figure}
    \centering
    \includegraphics[width=0.475\textwidth]{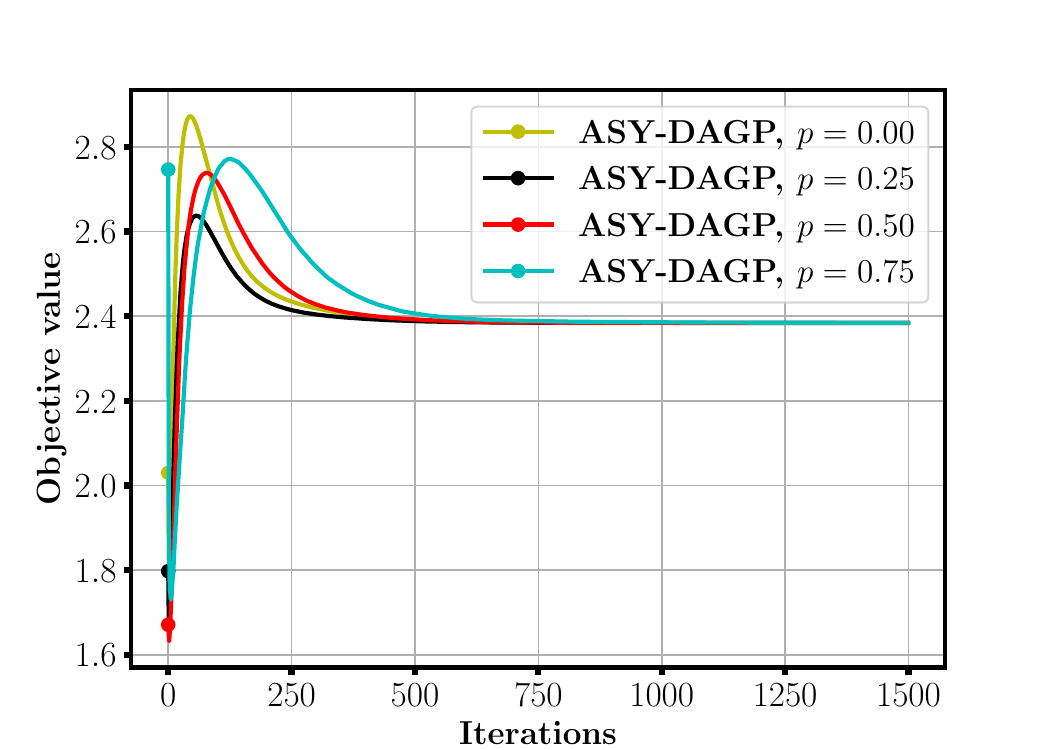}
    \caption{Robustness to message losses with the communication failure probability of $p$.}
    \label{fig: drop masg}
\end{figure}

\textbf{Communication failures:} The goal of this experiment is to demonstrate the effectiveness of ASY-DAGP in handling dropped messages.
We consider the same problem and parameters as in the first experiment. 
The messages can be lost with a failure probability of $p$. Figure~\ref{fig: drop masg} investigates  different values of $p$. We observe that  increasing the  probability  failure leads to more iterations for convergence. We performed similar experiments on APPG, but it failed to converge. 

The conducted experiments provide an evidence of ASY-DAGP's efficiency under diverse conditions. ASY-DAGP performs robustly even under extensive communication failures, tolerating a loss of more than 75\% of messages. Additionally, ASY-DAGP surpasses its synchronous counterpart in wall-clock convergence speed. Despite being originally designed for constrained setups with directed communication networks, it can be successfully used for either unconstrained problem scenarios or setups with undirected networks. 


\section{Related works}
\label{sec: litearture}
This review primarily focuses on decentralized gossip algorithms with general directed communication networks. 
Given the real-world relevance, complex dynamics, and asymmetric relationships inherent in directed graphs, they offer a more comprehensive and realistic model for decentralized optimization than undirected graphs.
%
Notable papers, such as \citep{lian2018asynchronous}, \citep{peng2016arock}, \citep{wu2023delay}, and \citep{even2021asynchronous}, have explored undirected communication graphs in asynchronous setups. However, their approach cannot be generalized to directed networks.
%
In the following,  we review synchronous optimization algorithms, followed by their asynchronous counterparts.

\textbf{Synchronous optimization algorithms:}  
Earlier papers consider the unconstrained version of \eqref{eq: p1}. For example,
\cite{nedic2014distributed, nedic2016stochastic} propose the \emph{subgradient-push} algorithm based on the push-sum protocol~\citep{kempe2003gossip}, which
requires a vanishing step size.
Taking a fixed step size leads to an error proportional to the average dissimilarity between local objective functions, a quantity related to  
\emph{distribution shift} in federated learning literature~\citep{reisizadeh2020robust, fallah2020personalized}. 
In response, a family of algorithms based on the so-called gradient tracking technique is introduced  \citep{xi2017dextra, nedic2017achieving}. In particular, the SONATA algorithm  proposed by \cite{scutari2019distributed} combines the push-sum protocol with a gradient tracking technique.  The Push-Pull algorithm proposed by \citet{xin2018linear} and \citet{pu2020push} forgoes the push-sum protocol, making it a  simpler algorithm to analyse and implement. Both SONATA and Push-Pull enjoy a linear convergence rate when the local objective functions are smooth and strongly convex. 
\citet{xi2016distributed} study constraints by proposing the DDPS algorithm which handles only identical constraints with diminishing step sizes.  The DAGP algorithm presented by \citet{firoozTAC} is the state-of-the-art, ensuring provable convergence for the generic case of \eqref{eq: p1} with a fixed step size and directed graphs.

\textbf{Asynchronous optimization algorithms:} 
The study by \citet{tsitsiklis1986distributed} is one of the first attempts to formulate the asynchronous setup, inspiring many other algorithms. See the comprehensive survey by 
\citet{assran2020survey}. 
\cite{zhang2019asyspa} and \citet{assran2020agp} propose the asynchronous version of the subgradient-push algorithm. They address the uneven update rate issue, which leads to biased solutions.
These methods require a decaying step size similar to their synchronous versions. 
The asynchronous version of the SONATA algorithm, ASY-SONATA, is proposed by \citet{tian2020asysonata}. They introduce the perturbed sum-push protocol with a gradient tracking technique.
Finally, \citet{zhang2019fully} propose the fully asynchronous push-pull gradient (APPG) algorithm.
Both ASY-SONATA and APPG achieve a linear rate of convergence using a constant step size: ASY-SONATA for strongly convex and smooth problems, while APPG for smooth problems satisfying the Polyak-Lojasiewicz condition.
%


Several papers address asynchrony and constraints by assuming composite objective functions with non-smooth terms, notably \citep{sayed2017decentralized} and \citep{latafat2022primal}. These studies are based on the dual formulation of the problem, leading to increased communication and computations due to the incorporation of dual variables. 
The algorithm by \citet{latafat2022primal} is not completely asynchronous; it only permits agents to handle delayed information.
\citet{sayed2017decentralized} consider undirected graphs,  and their algorithm requires step-size adjustments, which demands the knowledge of the agents' update rates — a stipulation not practical for real-world applications. Furthermore, it lacks a guaranteed rate of convergence.

\section{Conclusion}
We proposed an asynchronous modification of the DAGP algorithm that accommodates delays and missing messages, and is capable of solving smooth optimization problems with local constraints over directed graphs. The asynchronous updates eliminate idle time, thus leading to faster wall-clock convergence.
We presented a convergence analysis of our algorithm, based on a relaxation of the PEP framework, which we refer to as LQ-PEP. It is closely connected to the ALB approach in  DAGP analysis. Using LQ-PEP, we summarized the effect of the delay in a parameter that we referred to as delay response. Lastly, our experimental results substantiate the resilience against dropped messages, and show ASY-DAGP surpasses existing algorithms, even in their specialized, restricted setups.





\bibliography{refs}

\onecolumn
\newpage


\section*{{\huge{\hspace{4.5cm}Supplementary Material}}}
\vspace{0.1cm}

The equation numbers and references are in agreement with the main body of the paper.




\section{Mathematical Notation}\label{sec: mathematical_notation}
In this paper, we denote vectors using bold lowercase letters and matrices with bold uppercase letters. The element  at the $v^\tth$ row and $u^\tth$ column of matrix $\bW$ is denoted by $w_{v u}$. The transpose of $\bW$ is $\bW^{T}$, and its right null space is shown by $\ker(\bW)$. Therefore, $\bx$ belongs to $\ker(\bW)$ if and only if $\bW\bx = \bzero$.
We use $\bone_n$ and $\bzero_n$ to respectively denote $n-$dimensional vectors entirely consisting of ones and zeros. Furthermore, an $m \times n$ matrix containing only zero elements is shown by $\bO_{m\times n}$. The indices $m$ and $n$ may be omitted if there is no risk of confusion.
The Euclidean inner product of vectors is denoted by $\langle.,.\rangle$, and the matrix inner product is denoted by $\langle\bA,\bC\rangle = \tr(\bA\bC^T)$, where $\tr(\ldotp)$ signifies the matrix trace operator. The Kronecker delta function is shown by $\delta_{k,l}$. $\partial I_S$ indicates the normal cone of $S$. $|T|$ indicates the cardinality of set $T$. 

In this paper, we generally use subscripts to define the iteration number, and superscripts to denote the node number. For example, $\nabla f^v(\bx^v_{k})$ represents the gradient of the local function of node $v$ at its local variable at the $k^\tth$ iteration. Additionally, to provide matrix representations of equations, we arrange vector variables as the rows of a matrix. For instance, the matrix $\bG\in \bbR^{M\times m}$ contains all $\bg^v \in \bbR^m $ vectors, with $v\in\{1,\ldots,M\}$, as its rows. For the sake of simplicity, in this paper, the statement $\bG\in\ker(\bW)$ implies that each column of $\bG$ is an element in the null space of $\bW$.

\section{Analysis of fixed points}\label{sec: fixed-point}
The fixed point analysis is a simple way to study the convergence properties of optimization algorithms. In DAGP, the extra assumption $\ker(\bQ) = \ker(\bW^T)$ is required to show any fixed point of DAGP is a consensus and optimal solution of problem \eqref{eq: p1}. However, In ASY-DAGP, by defining new $\bp^v$ variables and designing their update equation, defined in \eqref{eq: asy-p_update}, based on the satisfaction of the optimality condition for the fixed point, there is no need for this extra assumption on gossip matrices anymore. For the fixed point iteration of ASY-DAGP, we have $\bx^v_{k+1} = \bx^v_{k} = \bx^v$. In addition, the same relation holds true for $\bg^v, \bp^v,$ and  $\bh^v$ variables. Since the algorithm is in the fixed point regime, we have additional relation $\ba^{vu}_{k+1} = \ba^{vu}_k = \bx^u$, and $\bb^{vu}_{k+1} = \bb^{vu}_k = \bp^u.$ Let us look at the fixed point iteration of ASY-DAGP in the following matrix form:
\begin{align}
    & \bZ = \bX-\bW\bX-\mu(\bnabla \bF-\bG) \label{eq: fix-z} \\
    & \bX = \bP(\bZ) \label{eq: fix-x}\\
    & \alpha (\bH-\bG) + \rho\left(\bnabla\bF-\bG+\frac{1}{\mu}(\bZ-\bX)\right) = \bO \label{eq: fix-g} \\
    & \bQ\bP = (\gamma-1)\bG \label{eq: fix-p} \\
    & \bQ\bP = (\gamma-1)\bH \label{eq: fix-h}
\end{align}

From the last two equations, $\bH=\bG$.  Then, from \eqref{eq: fix-g}:
\begin{equation} \label{eq: fix-g-2}
    \bG = \bnabla\bF+\frac{1}{\mu}(\bZ-\bX).
\end{equation}

This relation proves that the fixed point is the consensus and optimal solution of \eqref{eq: p1}. Using \eqref{eq: fix-g-2} in \eqref{eq: fix-z}, we have $\bW\bX = \bO$, which corroborate a consensus solution, i.e. $\bx^v=\bx^u=\bx$ for all $ v,u\in\calV$. To show  $\bx$ is an optimal solution, first note that $\bone^T\bG= \bzero^T$ from \eqref{eq: fix-p}. Moreover, consider that $\bz^v-\bx^v \in \partial I_{S^v}$, for all $v\in\calV.$
Then, by left multiplying \eqref{eq: fix-g-2} with $\bone^T$, considering the conic property of normal cone, we have the optimality condition satisfied at $\bx$.
By this analysis, we showed that ASY-DAGP does not require the limiting kernel assumption on gossip matrices. 

\newpage
\section{Convergence analysis}\label{sec: convergence}

In this section, we follow the steps of ALB, similar to~\citep{firoozTAC}, to derive LQ-PEP.

As mentioned in the paper, for the sake of analysis, we index the iterations of each node, individually. Hence, iteration $k$ at different nodes happen at different real times. However, we pretend that each iteration occurs simultaneously at every node, but the messages experience a corrected delay that may be negative (non-causal). This assumption does not affect the generality of our approach. 
Based on this convention, we introduce the local index sets $T^{vu}_k$ containing all iteration numbers $k^\prime$, where a message from $u$ is sent  at $k^\prime$ and is processed by $v$ at time $k+1$. In other words, $T^{vu}_k = \{k^\prime | k^\prime+\tau^{vu}_{k^\prime}=k+1\}$.

\subsection*{Step one: combining inequalities from the assumptions} 
First, we define
\begin{align*}
    F^v(\bx) &\coloneqq f^v(\bx)-f^v(\bx^*)-\langle\nabla f^v(\bx^*),\bx-\bx^*\rangle, \\
    T^v(\bx) &\coloneqq -\langle\bn^v,\bx-\bx^*\rangle,
\end{align*}
where $\bn^v \in \partial I_{S^v}(\bx^*)$, i.e. $\bn^v$ is an element in the normal cone of $S^v$ at $\bx^*$. From the convexity assumption on $f^v$ and the definition of normal cone, we have 
\begin{alignat}{2}
F^v(\bx^v) &\geq 0, && \qquad \forall \bx \\   
T^v(\bx^v) &\geq 0. && \qquad \forall \bx \in S^v   
\end{alignat}
Hence, these functions may serve as performance measures for our decentralized and constrained setup. 
In the following,  we define $F^v_{k+1}=F^v(\bx^v_{k+1})$ and  $T^v_{k+1}= T^v(\bx^v_{k+1})$ for notation simplicity.

Form convexity of $f^v$, we have
\begin{equation}\label{eq: 1}
    F^v_k+\left\langle\nabla f^v(\bx^*)-\nabla f^v(\bx^v_k), \bx^v_k-\bx^*\right\rangle \leq 0.
\end{equation}
From the $L-$smoothness property of $f^v$, we have 
\begin{equation}\label{eq: 2}
     F_{k+1}^v-F_k^v-\langle\nabla f^v(\bx_k^v)-\nabla f^v(\bx^*),\bx_{k+1}^v-\bx^v_k\rangle \leq 0. 
\end{equation}
From the projection step in \eqref{eq: asy-z_update}, and $\bx^v_{k+1}\in S^v$, we have
\begin{equation}\label{eq: 3}
    \mu T^v_{k+1}+\langle\bx^*-\bx^v_{k+1},\bz_{k+1}^v-\bx^v_{k+1}-\mu\bn^v\rangle\leq 0.
\end{equation}

By summing the first two inequalities, multiplying the resulting expression by $\mu$, and then adding it to \eqref{eq: 3}, we obtain
\begin{align}
    \mu(F^v_{k+1}+T^v_{k+1}) &-  \frac{L\mu}{2}\left\|\bx^v_{k+1}-\bx^v_k\right\| \nonumber \\
    &+ \Big\langle\bx^*-\bx^v_{k+1},\bz^v_{k+1}-\bx^v_{k+1} + \mu\big(\nabla f^v(\bx^v_k)-\nabla f^v(\bx^*) - \bn^v\big)\Big\rangle \leq 0.
    \label{eq: 4}
\end{align}

\subsection*{Step two: plugging algorithm's dynamics} 
The objective of this step is to eliminate the gradient term from \eqref{eq: 4}. By substituting the definition of $\bz^v_{k+1}$ from \eqref{eq: asy-z_update} into \eqref{eq: 4}, we obtain
\begin{align}
    \mu (F_{k+1}^v & +T_{k+1}^v)  -\frac {L\mu}2\left\|\bx_{k+1}^v-\bx_{k}^v\right\|^2  \nonumber  \\ 
    &+ \Big\langle\bx^*-\bx^v_{k+1}, 
     \bx_k^v-\suml_{u}w_{vu}\ba^{vu}_k-\bx^v_{k+1}  +\mu(\bg^v_k-\nabla f^v(\bx^*)-\bn^v)\Big\rangle \leq 0. \label{eq: sumall}
\end{align}

Furthermore, to remove the algorithm dynamics' dependency on $\bz^v_{k+1}$, we can plug the definition of $\bz^v_{k+1}$ into the update rule or the dynamics of $\bg^v_k$ as shown in \eqref{eq: g_update}. This gives us
\begin{equation}
    \bg_{k+1}^v=\bg_k^v+\frac\rho\mu\left(\bx_k^v-\suml_{u}w_{vu}\ba^{vu}_k-\bx^v_{k+1}\right)+\alpha(\bh^v_k-\bg^v_k).
\end{equation}

Now, for simplicity in analysis and notation, we  define $\tlbx_{k}^v\coloneqq \bx^v_k-\bx^*$, $\tlbg^v_k\coloneqq \bg^v_k-(\nabla f^v(\bx^*)+\bn^v)$, $\tlbh^v_k = \bh^v_k - (\nabla f^v(\bx^*) + \bn^v)$, and $\tlbp^v_k = \bp^v_k - \frac{\gamma-1}{\sum_uq_{vu}}(\nabla f^v(\bx^*) + \bn^v)$.
This transformation adjusts the system's origin to the optimal solution, leading to a modification of the equation in \eqref{eq: sumall} to 
\begin{equation}\label{eq: sumall_shift}
    \mu \left(F_{k+1}^v+T_{k+1}^v\right)  -\frac {L\mu}2\left\|\tlbx_{k+1}^v-\tlbx_{k}^v\right\|^2  -  \Big\langle \tlbx^v_{k+1}, 
     \tlbx_k^v-\suml_{u}w_{vu}\tlba^{vu}_k-\tlbx^v_{k+1} + \mu \tlbg\Big\rangle \leq 0. 
\end{equation}
Similarly, the introduction of new variables transforms the update rules for the algorithm's variables into
\begin{align}
\tlbg_{k+1}^v&=\tlbg_k^v+\frac\rho\mu\left(\tlbx_k^v-\suml_{u}w_{vu}\tlba^{vu}_k-\tlbx^v_{k+1}\right)+\alpha(\tlbh^v_k-\tlbg^v_k), \label{eq: tlg_update}\\
\tlbp^v_{k+1} &=  \tlbp^v_k - \eta\suml_{u} q_{vu}\tlbb^{vu}_k + \eta(\gamma-1)\tlbg^v_k, \label{eq: tlp_update}\\
\tlbh_{k+1}^v  &=  \gamma\tlbh_k^v -\suml_{u} q_{v u}\tlbb^{vu}_k, \label{eq: tlh_update}
\end{align}
\begin{equation} \label{eq: tla_update}
\tlba^{vu}_{k+1} = 
\left\{ 
\begin{array}{lc}
\tlba^{vu}_k &  \text{\footnotesize If the buffer is empty. 
}\\
\frac{1}{|T^{vu}_k|} \suml_{k^\prime\in T^{vu}_k} \tlbx^u_{k^\prime}  & \text{ \footnotesize Otherwise 
} 
\end{array}
\right.
\end{equation}
\begin{equation} \label{eq: tlb_update}
\tlbb^{vu}_{k+1} = 
\left\{ 
\begin{array}{lc}
\tlbb^{vu}_k &  \text{\footnotesize If the buffer is empty 
} \\
\frac{1}{|T^{vu}_k|} \suml_{k^\prime\in T^{vu}_k} \tlbp^u_{k^\prime}  & \text{ \footnotesize Otherwise 
} 
\end{array}
\right.
\end{equation}

By defining the state vectors for ASY-DAGP as $ \mathbf{\Psi}_k \in \mathbb{R}^{(2M^2+6M) \times m} $ 
\begin{equation}
\bPsi_k = \left[ \bX_{k+1}^T \;\; \bP_{k+1}^T \;\; \bX_k^T \;\; \bG^T_k \;\;  \bH^T_k \;\; \bP_k^T \;\; \{\bA_{v,k}^T\}_{v=1}^M \;\; \{ \bB_{v,k}^T \}_{v=1}^M  \right]^T,
\end{equation}
and introducing the shifted states $ \tilde{\mathbf{\Psi}}_k = \mathbf{\Psi}_k - \mathbf{\Psi}_* $, where $ \mathbf{\Psi}_* $ represents the optimal states, we can describe the evolution of ASY-DAGP's variables using the subsequent non-causal linear dynamical system
\begin{equation}\label{eq: non-causal-dynamic}
        \tlbPsi_{k+1} - \bar\bR\tlbPsi_k -  \bP\tilde\bU_{k} = \suml_{l=0}^{K-2} \tilde\bR_{k,k+1-l}\tlbPsi_l, \qquad k = 0,\ldots,K-2
\end{equation}
where $\barbR,\tlbR_{k,k+1-l}$, and  $\bP$ matrices are defined in \eqref{eq: barR_def}, \eqref{eq: tildeR_def}, and \eqref{eq: P_def}, respectively. The matrix $\tilde\bU_k \coloneqq [\tlbX_{k+2}^T \;\; \tlbP_{k+2}^T ]^T$ represents the input of the system. 

Regarding the creation of $\barbR$ and $\tlbR_{k,0}$ matrices, we follow the approach that the ideal case, where there are no communication delays and the agents have instant access to the variables of other nodes,  should be encapsulated within the definition of $\barbR$. By doing so, in scenarios without delays, the ASY-DAGP analysis simplifies to the DAGP analysis.
Therefore, we put all the delay-independent terms in  \eqref{eq: tlg_update}, \eqref{eq: tlp_update}, and \eqref{eq: tlh_update} in the definition of $\barbR$.
Also, in the creation of $\barbR$ the blocks representing update rules \eqref{eq: tla_update} and \eqref{eq: tlb_update} are created based on the ideal case, i.e. $\tlba^{vu}_{k+1}= \tlbx^{u}_{k+1}$ and $\tlbb^{vu}_{k+1}=\tlbp^u_{k+1}$.
Since the coefficient of $\tlbPsi_k$ is $\barbR + \tlbR_{k,0}$, and we put the ideal case in $\barbR$,  the $\tlbR_{k,0}$ matrix should compensate for non-ideal cases. 
In $\tlbR_{k,k+1-l}$, only the blocks associated with $\tlbA_{v,k}$ and $\tlbB_{v,k}$ are non-zero. Given the similarity of update rules in \eqref{eq: tla_update} and \eqref{eq: tlb_update}, we focus on how we create the blocks representing the update rule of $\tlbA_{v,k}$. For clarity, we consider only one node $v$ and the $u\tth$ row in $\tlbA_{v,k}$. For simplicity in notation, we define the \emph{interaction element} as the possible non-zero element connecting $\tlba^{vu}_{k+1}$ and $\tlbx^u_{k+1}$, which is in the $u^\tth$ column of the matrix $\tlbR_{k,k+1-l}$. 
When computing $\tlbR_{k,0}$, several scenarios arise:
\begin{itemize}
    \item If the buffer is empty and the ideal message is not present, interaction element is $-1$. Additionally, the element linking $\ba^{vu}_{k+1}$ and $\ba^{vu}_k$ is $+1$. 
    \item If several messages exist in the buffer and $k+1 \in T^{vu}_k$, i.e. ideal message is present in the buffer, the interaction element computes as $-1+\frac{1}{|T^{vu}_k|}$.
    \item If several messages are in the buffer, but $k+1$ is not in the indices set $T^{vu}_k$, the interaction element is $-1$.
\end{itemize}
For the computation of other $\tlbR_{k,k+1-l},\; l\neq k+1$, the only case when an element can be non-zero is when there are multiple messages in the buffer, excluding the ideal message. In this scenario, the element connecting $\tlba^{vu}_{k+1}$ and $\tlbx^u_{k+1}$ is $\frac{1}{|T^{vu}_k|}$ for all $l\in T^{vu}_k$.

\subsection*{Step three: defining performance measures}  
As mentioned earlier, $F^v$ and $T^v$ are potential performance measures for the decentralized constrained setup. However, in decentralized setups, convergence to a consensus solution is of interest as well. Therefore, we add and remove $\frac{\zeta}2\sum_{u,v}\|\tlbx_{k+1}^u-\tlbx_{k+1}^v\|^2 $ to \eqref{eq: sumall_shift}. Then, by summing over all $v=1,\ldots,M$ and $k=0,\ldots,K-1$, we have
\begin{align}\label{eq: sumall_final}
    & \suml_{k=0}^{K-1} \left( \suml_{v=1}^M\mu\left(F_{k+1}^v+T_{k+1}^v\right) \nonumber
    +\frac{\zeta}2\suml_{u,v=1}^M\|\tlbx_{k+1}^u-\tlbx_{k+1}^v\|^2 \right)  \nonumber\\
    & -\frac{L\mu}2 \suml_{k=0}^{K-1} \suml_{v=1}^M \left\|\tlbx_{k+1}^v-\tlbx_{k}^v\right\|^2   \nonumber\\ 
    & - \suml_{k=0}^{K-1}\suml_{v=1}^M \Big\langle \tlbx^v_{k+1}, \tlbx_k^v-\suml_{u}w_{vu}\tlba^{vu}_k-\tlbx^v_{k+1} + \mu \tlbg\Big\rangle \nonumber\\
    & -\frac{\zeta}2 \suml_{k=0}^{K-1} \suml_{u,v=1}^M\|\tlbx_{k+1}^u-\tlbx_{k+1}^v\|^2   \leq 0.
\end{align}
The first line in \eqref{eq: sumall_final} includes all terms that measure the algorithm's performance over all iterations, represented by $\sum_k\Phi(\bPsi_k, \bPsi_*)$. In the Performance Estimation Problem (PEP), our goal is to find the worst-case performance. We define all terms in the last three summations in \eqref{eq: sumall_final} as $\bar A_K$, which can be written in the following quadratic form
\begin{equation}
    \bar A_K = \sum_{k=0}^{K-1} \langle \tlbPsi_k,\bS\tlbPsi_k \rangle,
\end{equation}
where matrix $\bS$ is defined in \eqref{eq: S_define}.  
Now, by applying the new definitions into \eqref{eq: sumall_final}, we have 
\begin{equation}
    \suml_{k=0}^{K-1} \Phi(\bPsi_k,\bPsi_*) \leq - \bar A_K.
\end{equation}
Hence, $\bar A_K$ is an upper bound for the performance measure, and the objective function in PEP can be replaced by $-\bar A_K$.
We then substitute the problem of maximizing $\Phi$ with the problem of minimizing $\bar A_K$, which has a quadratic form. We also consider the constraint as the linear dynamical system in \eqref{eq: non-causal-dynamic}, which contain all possible states created by the optimization algorithm. 
We call the new problem LQ-PEP. 
Upon reviewing the derivation of \eqref{eq: sumall_final}, we see that it uses a set of $\calF-$interpolation inequalities introduced by \citet  {taylor2017smooth}, representing the set of convex and smooth functions. Therefore, LQ-PEP is just a relaxation of PEP.

\subsection*{Step four: solving LQ-PEP} To establish that $\bar A_K$ is lower bounded by $-C\nf{\tlbPsi_0}^2$, which is essentially the solution to LQ-PEP, we need to show that zero is the solution to the following optimization problem.
\begin{eqnarray}
    &\minl_{\{\tilde\bPsi_k\}_{k=0}^{K-1},\{\tilde\bU_k\}_{k=0}^{K-2}} \;\;
    \suml_{k=0}^{K-1} \langle \tilde\bPsi_k,\bS\tilde\bPsi_k\rangle+C\|\tilde\bPsi_0\|^2\nwl
    &\text{s.t.}\nwl
    & \tlbPsi_{k+1} - \bar\bR\tlbPsi_k -  \bP\tilde\bU_{k} = \suml_{l=0}^{K-2} \tilde\bR_{k,k+1-l}\tlbPsi_l. \qquad k = 0,\ldots,K-2
    \label{eq: ALB_opt_again}
\end{eqnarray}
In the reference \citep{firoozTAC}, the dynamics of the algorithm are depicted as a linear time-invariant (LTI) system, as there are no communication delays. This distinction increases the complexity of the analysis, making our study substantially different from their analysis.

As this step is rather substantial and possibly intriguing to readers on its own, we will apply the result $\bar A_K \geq -C\nf{\tlbPsi_0}^2$ for now and provide a full proof after the convergence rates are determined.

\subsection*{Step five: providing the rates of convergence} Since $\bar A_K$ has a lower-bound $-CC_0$, where $C_0\coloneqq \nf{\tlbPsi_0}^2$, we obtain the following rates of convergence. 

     \textit{\textbf{Consensus:}} 
    Time-averaged local solutions converge to the consensus solution $\barbx_K$ as 
    \begin{equation}
            \|\barbx_K-\barbx^v_K\|_2^2\leq\frac{C_0C}{ M K}.  \qquad \forall v \in \calV
    \end{equation}
    
     \textit{\textbf{Feasibility gap:}}
    The consensus solution $\barbx_K$ goes into each constraint set with the following rate. 
    Hence, the distance between $\barbx_K$ and the feasible set decays at the same rate.
    \begin{equation}
    \mathrm{dist}^2(\barbx_K, S^v)=\calO(\frac 1K) \qquad \forall v \in \calV 
    \end{equation}
    
     \textit{\textbf{Optimality gap:}}
    The objective value converges to the optimal value as
    \begin{equation}
    \left|\suml_v f^v(\barbx^v_K)-\suml_v f^v(\bx^*)\right|\leq 
    \frac{C_0C}{\mu K}+\sqrt{\frac{C_0CC_2}{ M K}},
    \end{equation} 
    where  $C_2= \sqrt{\sum_v\|\bn^v+\nabla f^v(\bx^*)\|^2}$, and it only depends on the choice of the regular optimal solution (which exists by the assumptions).
\qed

\newpage
\vspace*{\fill}
\begin{center}
{\Large{\textbf{Matrix definitions:}}}
\end{center}
We define matrix $\bar\bW \coloneqq \text{blkdiag}(\bW)$ as a block-diagonal matrix where each block on its diagonal is given by a single row from $\bW.$ Similarly, we define matrix $\bar\bQ \coloneqq \text{blkdiag}(\bQ)$. 

\begin{center}
    \begin{equation}  \label{eq: barR_def}
    \barbR = \left[  
    \begin{array}{cccccccc}
         \bO& 
         \bO& 
         \bO&
         \bO&
         \bO&
         \bO&
         \bO&
         \bO
         \\
         \bO& 
         \bO& 
         \bO&
         \bO&
         \bO&
         \bO&
         \bO&
         \bO
         \\
         \bI&
         \bO& 
         \bO&
         \bO&
         \bO&
         \bO&
         \bO&
         \bO
         \\
         -\frac{\rho}{\mu}\bI&
         \bO& 
         \frac{\rho}{\mu}\bI&
         (1-\alpha)\bI&
         \alpha\bI&
         \bO&
         -\frac{\rho}{\mu}\bar\bW &
         \bO
         \\
         \bO&
         \bO& 
         \bO&
         \bO&
         \gamma\bI&
         \bO&
         \bO&
         \bar\bQ
         \\
         \bO&
         \bO& 
         \bO&
         \eta(\gamma-1)\bI&
         \bO&
         \bI&
         \bO&
         -\eta
         \bar\bQ
          \\
         \bone_M\otimes\bI &
         \bO& 
         \bO_{M^2\times M}  &
         \bO&
         \bO&
         \bO&
         \bO&
         \bO
          \\
         \bO&
         \bone_M\otimes\bI &
         \bO&
         \bO&
         \bO&
         \bO&
         \bO&
         \bO
    \end{array}
    \right]
\end{equation}

\begin{equation}
    \label{eq: tildeR_def}
    \left\{\tlbR_{k,k+1-l}\right\}_{i,j} = 
    \left\{ 
    \begin{array}{cl}
    -1 & \quad   l =k+1, |T^{vu}_k| = 0, i=(v+5)M+u, j=u\\ 
    +1 & \quad   l = k+1, |T^{vu}_k| = 0, i=(v+5)M+u, j = (v+5)M+ u\\ 
    -1+\frac{1}{|T^{vu}_k|} & \quad l = k+1,  |T^{vu}_k| \neq 0, l\in T^{vu}_k, i=(v+5)M+u, j=u \\ 
    -1 & \quad  l =k+1,  |T^{vu}_k| \neq 0,  l \notin T^{vu}_k, i=(v+5)M+u, j=u \\ 
    \frac{1}{|T^{vu}_k|} & \quad l\neq k+1, |T^{vu}_k| \neq 0, l \in T^{vu}_k, i=(v+5)M+u, j=u \\
    ---- & \quad ------ \\
    -1 & \quad   l =k+1, |T^{vu}_k| = 0, i=M^2+(v+5)M+u, j=M+u\\ 
    +1 & \quad   l = k+1, |T^{vu}_k| = 0, i=M^2+(v+5)M+u, j = M^2+(v+6)M+ u\\ 
    -1+\frac{1}{|T^{vu}_k|} & \quad l = k+1,  |T^{vu}_k| \neq 0, l\in T^{vu}_k, i=M^2+(v+5)M+u, j=M+u \\ 
    -1 & \quad  l =k+1,  |T^{vu}_k| \neq 0,  l \notin T^{vu}_k, i=M^2+(v+5)M+u, j=M+u \\ 
    \frac{1}{|T^{vu}_k|} & \quad l\neq k+1, |T^{vu}_k| \neq 0, l \in T^{vu}_k, i=M^2+(v+5)M+u, j=M+u \\
    ---- & \quad ------ \\
     0 & \quad   \text{otherwise}
    \end{array}
    \right.
\end{equation}

\begin{equation}\label{eq: P_def}
    \bP = \left[  
    \begin{array}{cc}
         \bI_{M\times M} &   \bO_{M\times M}  \\
         \bO_{M\times M} &   \bI_{M\times M} \\ 
         \bO_{(2M^2+4M)\times M}  & \bO_{(2M^2+4M)\times M}
    \end{array}
    \right]
\end{equation}

\begin{equation}\label{eq: S_define}
    \bS=\left[
    \begin{array}{cccccccc}
        \left(1-\frac{L\mu}2\right)\bI-M\zeta\left(\bI-\frac 1M\bone\bone^T\right)& 
        \bO &
        -\frac {L\mu-1}{2}\bI & 
        -\frac {\mu}2\bI  & 
        \bO &
        \bO &
        \frac 12
        \bar\bW &
        \bO 
        \\
        \bO &
        \bO &
        \bO &
        \bO &
        \bO &
        \bO &
        \bO &
        \bO
        \\ 
        -\frac {L\mu-1}{2}\bI&
        \bO &
        -\frac{L\mu}2\bI &
        \bO &
        \bO &
        \bO&
        \bO&
        \bO 
        \\
        -\frac{\mu}2\bI &
        \bO &
        \bO & 
        \bO & 
        \bO &
        \bO&
        \bO&
        \bO 
        \\
        \bO &
        \bO & 
        \bO & 
        \bO &
        \bO&
        \bO&
        \bO &
        \bO
        \\
        \bO &
        \bO & 
        \bO & 
        \bO &
        \bO&
        \bO&
        \bO &
        \bO
        \\
        \frac 12\bar\bW &
        \bO &
        \bO & 
        \bO & 
        \bO &
        \bO &
        \bO &
        \bO 
        \\
        \bO &
        \bO & 
        \bO & 
        \bO &
        \bO&
        \bO&
        \bO &
        \bO
        \end{array}
    \right]
\end{equation}

\end{center}
\vspace*{\fill}


\newpage
\section{Solving LQ-PEP }\label{sec: lqpep}
In this proof, $D$ denotes a generic constant whose value might increase each time it appears.

In order to solve the LQ-PEP and demonstrate that $-C\nf{\tlbPsi_0}^2 \leq \bar A_K$, we need to establish that zero is the optimal value of \eqref{eq: ALB_opt_again}. Notably, $\bS$ being indefinite sets LQ-PEP apart from linear quadratic control. If zero is not the optimum of that optimization, the problem becomes unbounded, pushing the optimal value to $-\infty$. Consequently, we restrict the optimization problem in \eqref{eq: ALB_opt_again} to a unit circle. This constraint guides us towards the following optimization:
\begin{eqnarray}
    &\minl_{\{\tilde\bPsi_k\}_{k=0}^{K-1},\{\tilde\bU_k\}_{k=0}^{K-2}}
    \suml_{k=0}^{K-1} \frac{1}{2}\langle \tilde\bPsi_k,\bS\tilde\bPsi_k\rangle+\frac{C}{2}\|\tilde\bPsi_0\|^2\nwl
    &\text{s.t.}\nwl
    & \tlbPsi_{k+1} - \bar\bR\tlbPsi_k -  \bP\tilde\bU_{k} = \suml_{l=0}^{K-2} \tilde\bR_{k,k+1-l}\tlbPsi_l, \qquad k = 0,\ldots,K-2
    \label{eq: ALB_opt_unitcircle} \nwl
    &\frac{1}{2}\nf{\tlbPsi_0}^2 + \frac{1}{2}\suml_{k=0}^{K-2}\|\tilde\bU_{k}\|_\mathrm{F}^2 \leq \frac{1}{2}.\qquad\qquad   \label{eq: alb_constrained}
\end{eqnarray}
If the optimal value for this new optimization is negative, the original optimization in \eqref{eq: ALB_opt_again} becomes unbounded. Thus, asserting that zero is the optimal solution for the original problem is equivalent to stating that the optimal value of the optimization in \eqref{eq: alb_constrained} is non-negative.

To show that the optimal value of \eqref{eq: alb_constrained} is non-negative,  we introduce the dual Lagrangian multipliers $\bLambda_k$ and $\beta\geq0$. Consequently, we construct the dual optimization problem. Upon solving this problem with respect to both the primal and dual variables, we achieve
\begin{alignat}{2}
    \tlbPsi_{k+1} - \bar\bR\tlbPsi_k - \bP\tilde\bU_k &= \suml_{l=0}^{K-2} \tilde\bR_{k,k+1-l}\tlbPsi_l &&\qquad k = 0,\ldots,K-2 \label{eq: recurrence_eqs_1} \\
    \bLambda_{k-1} - \bar\bR^T\bLambda_k + \bS\tlbPsi_k &= \suml_{l=0}^{K-2} \tilde\bR^T_{l,l-k+1}\bLambda_l &&\qquad k = 0,\ldots,K-2
    \label{eq: recurrence_eqs_2} \\
    \beta\tilde\bU_k - \bP^T\bLambda_k & = \bO && \qquad k=0,\ldots,K-2 \label{eq: recurrence_eqs_3}
\end{alignat}
with boundary conditions 
\begin{align}
    \bLambda_{-1} + (\beta+C)\tlbPsi_0 & = \bO  \label{eq: boundaries_1} \\
    \bLambda_{K-2} + \bS\tlbPsi_{K-1}  & = \bO  \label{eq: boundaries_2}
\end{align}
After some simplifications and substitutions, we notice that the optimal value of this problem $-\beta\left( \nf{\bPsi_0}^2 + \sum_{k=0}^{K-2}\nf{\tilde\bU_k}^2\right)$. Since $\beta$ is positive, the optimal solution is negative unless the system of linear recurrences  in \eqref{eq: recurrence_eqs_1}, \eqref{eq: recurrence_eqs_2}, and \eqref{eq: recurrence_eqs_3} has no non-zero solution. 


In summary,  if we show the system of linear recurrences  in \eqref{eq: recurrence_eqs_1}, \eqref{eq: recurrence_eqs_2}, and \eqref{eq: recurrence_eqs_3} has no non-zero solution, zero will be the optimal value of the optimization problems in \eqref{eq: ALB_opt_again} and \eqref{eq: alb_constrained}, consequently,  $\bar A_K$ has a lower bound $-C\nf{\tlbPsi_0}^2.$



To show this, first, we remove the dependence on $\tilde\bU_k$. Then, we have the new system of linear recurrences
\begin{alignat}{2}
    \tlbPsi_{k+1} - \bar\bR\tlbPsi_k - \frac{1}{\beta} \bP\bP^T\bLambda_k &= \suml_{l=0}^{K-2} \tilde\bR_{k,k+1-l}\tlbPsi_l &&\qquad k = 0,\ldots,K-2 \label{eq: recurrence_eqs_11} \\
    \bLambda_{k-1} - \bar\bR^T\bLambda_k + \bS\tlbPsi_k &= \suml_{l=0}^{K-2} \tilde\bR^T_{l,l-k+1}\bLambda_l &&\qquad k = 0,\ldots,K-2
    \label{eq: recurrence_eqs_22}
\end{alignat}
with the same boundary conditions as in \eqref{eq: boundaries_1} and \eqref{eq: boundaries_2}. 
By defining the finite duration $z-$transforms of $\tlbPsi_k$, and $\bLambda_k$ as
\begin{equation}
    \tlbPsi(z) = \suml_{k=0}^{K-1} \tlbPsi_k z^k, \qquad \bLambda(z) = \suml_{k=0}^{K-1} \bLambda_{k-1}z^k,
\end{equation}
and taking the $z-$transform of \eqref{eq: recurrence_eqs_11} and \eqref{eq: recurrence_eqs_22}, we have
\begin{equation} 
\label{eq: 7}
\small
    \left[ 
    \begin{array}{c}
        \tlbPsi(z)\\
        \bLambda(z)
    \end{array}
    \right] = \bF_\beta^{-1}(z) \left\{ \left[ 
    \begin{array}{c}
         \bS\tlbPsi_{K-1}z^{K-1} + \bLambda_{K-2}z^{K-1} - z^{-1}\bar\bR^T\bLambda_{-1} \\
          \tlbPsi_0 - \bar\bR\tlbPsi_{K-1}z^K - \frac{1}{\beta}\bP\bP^T\bLambda_{-1}
    \end{array}
    \right]
    + 
    \left[
    \begin{array}{c}
         \suml_{k,l=0}^{K-2} \tilde\bR^T_{l,l-k+1}\bLambda_l z^k  \\
          \suml_{k,l=0}^{K-2} \tilde\bR_{k,k+1-l}\tlbPsi_l z^{k+1}
    \end{array}
    \right]
    \right\}
\end{equation}
where
\begin{equation}
    \bF_\beta(z) = \left[
    \begin{array}{cc}
         \bS & \bI-z^{-1}\bar\bR^T \\
         \bI-z\bar\bR & -\frac{1}{\beta}\bP\bP^T
    \end{array}
    \right].
\end{equation}

From the definition of inverse $z-$transform, we have
\begin{align}
    \tlbPsi_m &= \frac{1}{2\pi j} \; \oint z^{-(m+1)}\tlbPsi(z)\mathrm{d}z, \\
    \bLambda_{m-1} &= \frac{1}{2\pi j} \; \oint z^{-(m+1)} \bLambda(z) \mathrm{d}z.
\end{align}
%
%
Given the structure of $\bF_\beta(z)^{-1}$ from Lemmas \ref{lemma: matrix_decomposition} and \ref{lemma: inversly-related_roots}, we compute the inverse z-transforms, based on Lemma \ref{lem:integral}, as:
\begin{equation}\label{eq: solution_recurrence}
    \left[  
    \begin{array}{c}
         \tlbPsi_m  \\
         \bLambda_{m-1} 
    \end{array}
    \right] \; = \;
    \suml_{p=1}^q \ba_p\bb_p^T
    \left[ 
    \begin{array}{c}
         \big(  \bS\tlbPsi_{K-1}+\bLambda_{K-2}\big)z_p^{K-m-2} + \barbR^T\bLambda_{-1}z_p^{m+2} \\
         \big(\frac{1}{\beta} \bP\bP^T\bLambda_{-1}-\tlbPsi_0\big)z_p^{m+1} - \barbR\tlbPsi_{K-1}z_p^{K-m-1}
    \end{array}
    \right] \; +  \;     
    \left[
    \begin{array}{c}
         \bar\bPsi_m  \\
         \bar\bLambda_{m-1} 
    \end{array}
    \right],
\end{equation}
with
%
%
\begin{equation}\label{eq: 57}
    \left[
    \begin{array}{c}
         \bar\bPsi_m  \\
         \bar\bLambda_{m-1} 
    \end{array}
    \right] \; = \;
    \suml_{p=1}^q \ba_p\bb_p^T
    \left[ 
    \begin{array}{c}
          \suml_{l=0}^{K-2} \left( \suml_{k=m+1}^{K-2} \tilde\bR^T_{l,l-k+1} z_p^{|k-m-1|}  
         - \suml_{k=0}^{m} \tilde\bR^T_{l,l-k+1} z_p^{|k-m-1|}\right)  \bLambda_l \\
          \suml_{l=0}^{K-2} \left( \suml_{k=m}^{K-2} \tilde\bR_{k,k+1-l}  z_p^{|k-m|} 
         - \suml_{k=0}^{m-1} \tilde\bR_{k,k+1-l} z_p^{|k-m|} \right)\tlbPsi_l
    \end{array}
    \right], 
\end{equation}
where $\ba_p \in \ker (\bF_\beta(z_p))$ and is assumed to be in the form of $ \ba_p = \left[ \bn_{p,\tlbPsi}^T \;\; \bn_{p,\bLambda}^T \right]^T$. 
The set $\{z_1,\ldots,z_q\}$ contains the poles of  $\bF_\beta(z)^{-1}$, which are inside the unit circle. 
We interpret the result in \eqref{eq: solution_recurrence} as the solution to the linear system of recurrence equations represented in \eqref{eq: recurrence_eqs_11} and \eqref{eq: recurrence_eqs_22}. This solution consists of two parts. The first part includes the solution to the homogeneous system, which is achieved by setting the right-hand side of equations \eqref{eq: recurrence_eqs_11} and \eqref{eq: recurrence_eqs_22} to zero. 
The second part, given in \eqref{eq: 57}, provides a particular solution, which is denoted by $\bar\bPsi_m, \bar\bLambda_{m-1}$ to simplify our further analysis.

For the sake analysis, we also write the homogeneous solution as
\[
\suml_{p=1}^q \ba_p\bgamma_p^Tz_p^{-m}\;,
\]
where
\begin{equation}
    \bgamma_p^T  =  \bb_p^T
    \left[ 
    \begin{array}{c}
         \big(  \bS\tlbPsi_{K-1}+\bLambda_{K-2}\big)z_p^{K-2} + \barbR^T\bLambda_{-1}z_p^{2m+2} \\
         \big(\frac{1}{\beta} \bP\bP^T\bLambda_{-1}-\tlbPsi_0\big)z_p^{2m+1} - \barbR\tlbPsi_{K-1}z_p^{K-1}
    \end{array}
    \right].
\end{equation}
This general solution can also be written as
\begin{equation} \label{eq: 56}
    \suml_{p=1}^q \ba_p\bgamma_p^Tz_p^{-m} = \left[  
    \ba_1z_1^{-m} \;\; \ldots \;\; \ba_qz_q^{-m} 
    \right] \left[ 
    \begin{array}{c}
         \bgamma_1^T  \\
         \vdots \\
         \bgamma_q^T
    \end{array}
    \right].
\end{equation}

The solution to the linear recurrence system must meet its boundary conditions. To achieve this, we substitute the solution into equations \eqref{eq: boundaries_1} and \eqref{eq: boundaries_2}. This substitution yields 
\begin{equation}
    \label{eq: apply_boundaries}
    \underbrace{\left[ 
    \begin{array}{cccc}
         \bn_{1,\bLambda} + (\beta+C)\bn_{1,\tlbPsi} & \ldots & \bn_{q,\bLambda} + (\beta+C)\bn_{q,\tlbPsi} \\
         (\bn_{1,\bLambda}+\bS\bn_{1,\tlbPsi})z_1^{-(K-1)} & \ldots & (\bn_{q,\bLambda}+\bS\bn_{q,\tlbPsi})z_q^{-(K-1)}
    \end{array}
    \right]}_{\bM} \left[
    \begin{array}{c}
        \bgamma_1^T  \\
        \vdots \\
        \bgamma_q^T
    \end{array}
    \right] = - \left[
    \begin{array}{c}
         \bar\bLambda_{-1} + (\beta+C)\bar\bPsi_0  \\
          \bar\bLambda_{K-2} + \bS\bar\bPsi_{K-1}
    \end{array}
    \right].
\end{equation}
To simplify, we define $\bM$ as illustrated in \eqref{eq: apply_boundaries}. Hence, for the solution to meet the boundary conditions, the following relation should be satisfied
\begin{equation} \label{eq: 55}
    \left[
    \begin{array}{c}
        \bgamma_1^T  \\
        \vdots \\
        \bgamma_q^T
    \end{array}
    \right] = -\bM^{-1} 
    \left[
    \begin{array}{c}
         \bar\bLambda_{-1} + (\beta+C)\bar\bPsi_0  \\
          \bar\bLambda_{K-2} + \bS\bar\bPsi_{K-1}
    \end{array}
    \right].
\end{equation}

Moving forward, the last step is to impose bounds on the aforementioned equalities to illustrate that these equations cannot be concurrently true except when zero is the optimal solution. 
To accomplish this, we make the equations simpler by looking at them column by column. Without loss of generality, If we consider any arbitrary column of $\tlbPsi$ and $\bLambda$, , and demonstrate that this column equals zero, then all the columns would be zero, completing our proof. For this reason, by $\norm{\tlbPsi}$, we are referring to a proper vector norm associated to one random column of $\tlbPsi$. Using this terminology and considering the $ \ell^\infty $ vector norm, we define
\begin{align}
    \lambda & = \; \max_{m=0,\ldots,K-1} \; \max  \left\{ 
    \| \tlbPsi_m \|_\infty,  \| \bLambda_{m-1}\|_\infty
    \right\}, \label{eq: lamda}\\
    \bar\lambda & = \; \max_{m=0,\ldots,K-1} \; \max  \Big\{ 
    \| \bar\bPsi_m \|_\infty,  \| \bar\bLambda_{m-1}\|_\infty
    \Big\}. \label{eq: lamdabar}
\end{align}

Putting all equations in \eqref{eq: solution_recurrence}, \eqref{eq: 57}, \eqref{eq: 56}, and \eqref{eq: 55} together, we get
\begin{equation}\label{eq: equation1_to_bound}
    \left[  
    \begin{array}{c}
         \tlbPsi_m  \\
         \bLambda_{m-1} 
    \end{array}
    \right] \; = \;
    -\left[  
    \ba_1z_1^{-m} \;\; \ldots \;\; \ba_qz_q^{-m} 
    \right] \bM^{-1} 
    \left[
    \begin{array}{c}
         \bar\bLambda_{-1} + (\beta+C)\bar\bPsi_0  \\
          \bar\bLambda_{K-2} + \bS\bar\bPsi_{K-1}
    \end{array}
    \right] 
    \; + \; 
     \left[
    \begin{array}{c}
         \bar\bPsi_m  \\
         \bar\bLambda_{m-1} 
    \end{array}
    \right].
\end{equation}
Taking norm, considering \eqref{eq: lamda} and \eqref{eq: lamdabar}, and using the triangle inequality, we get
\begin{equation} \label{eq: 76}
    \lambda  \leq   \left\|    
    \left[  
    \ba_1z_1^{-m} \;\; \ldots \;\; \ba_qz_q^{-m} 
    \right] \bM^{-1} 
    \left[
    \begin{array}{c}
         \bar\bLambda_{-1} + (\beta+C)\bar\bPsi_0  \\
          \bar\bLambda_{K-2} + \bS\bar\bPsi_{K-1}
    \end{array}
    \right]
    \right\| \; + \;  \bar\lambda.
\end{equation}

From \eqref{eq: apply_boundaries} and \eqref{eq: 55}, one can show that for sufficiently large $K$, there exists a generic constant $D$, dependant on constant $C$,  such that $\| \bgamma_p^T \|_\infty \leq D|z_p|^K \bar\lambda$, for all $p=1,\ldots,q$. 
%
Then, 
we have 
\begin{equation} \label{eq: 76_new}
    \lambda  \; \leq \;  \bar\lambda + D \bar\lambda \, \max_{m}\, \suml_{p=1}^q |z_p|^{K-m} \|\ba_p\| \;  \leq \;  D \bar\lambda 
\end{equation}
The above inequalities result in $\lambda \leq D \bar\lambda$. Now, we aim to determine an upper bound for $\bar\lambda$. Recall the definition in \eqref{eq: 57}
\begin{equation}\label{eq: 57_new}
    \left[
    \begin{array}{c}
         \bar\bPsi_m  \\
         \bar\bLambda_{m-1} 
    \end{array}
    \right] \; = \;
    \suml_{p=1}^q \ba_p\bb_p^T
    \left[ 
    \begin{array}{c}
          \suml_{l=0}^{K-2} \left( \suml_{k=m+1}^{K-2} \tilde\bR^T_{l,l-k+1} z_p^{|k-m-1|}  
         - \suml_{k=0}^{m} \tilde\bR^T_{l,l-k+1} z_p^{|k-m-1|}\right)  \bLambda_l \\
          \suml_{l=0}^{K-2} \left( \suml_{k=m}^{K-2} \tilde\bR_{k,k+1-l}  z_p^{|k-m|} 
         - \suml_{k=0}^{m-1} \tilde\bR_{k,k+1-l} z_p^{|k-m|} \right)\tlbPsi_l
    \end{array}
    \right], 
\end{equation}
Taking norm, considering \eqref{eq: lamda} and \eqref{eq: lamdabar},  we have 
%
\begin{align}
    \bar\lambda &  \leq   \max_m \suml_p \norm{\ba_p}_\infty \left| 
    \bb_p^T
    \left[ 
    \begin{array}{c}
          \suml_{l=0}^{K-2} \left( \suml_{k=m+1}^{K-2} \tilde\bR^T_{l,l-k+1} z_p^{|k-m-1|}  
         - \suml_{k=0}^{m} \tilde\bR^T_{l,l-k+1} z_p^{|k-m-1|}\right)  \bLambda_l \\
          \suml_{l=0}^{K-2} \left( \suml_{k=m}^{K-2} \tilde\bR_{k,k+1-l}  z_p^{|k-m|} 
         - \suml_{k=0}^{m-1} \tilde\bR_{k,k+1-l} z_p^{|k-m|} \right)\tlbPsi_l
    \end{array}
    \right]
    \right| 
    \label{eq: 811}
    \\
    & \leq D\lambda  \; \max_m \; \suml_p \; \max_l \;\; \max\{ \|\bA_{l,p,m} \|, \|\bB_{l,p,m} \| \} \label{eq: 83}
\end{align}
where
\begin{align}
    \bA_{k,p,m} & = \sum_{l=m+1}^{K-2} \tlbR^T_{k, k+1-l}z_p^{|l-m-1|} - \suml_{l=0}^m \tlbR^T_{k, k+1-l}z_p^{|l-m-1|} \\
    \bB_{l,p,m} & = b\sum_{k+m}^{K-2} \tlbR_{k,k+1-l}z_p^{k-m} - \suml_{k=0}^{m-1}\tlbR_{k,k+1-l}z_p^{k-m}
\end{align} 

where $\norm{A}_{1,\infty}$ is the $ \ell_1 - \ell_\infty $ mixed norm of matrix $ A $, given by
$ ||A||_{1,\infty} = \max_i \sum_j |a_{ij}|$.
To derive \eqref{eq: 83}, we applied the triangle inequality twice.
For \eqref{eq: 811}, we utilized the relationship 
$ | \mathbf{x}^T \mathbf{y} | \leq \|\mathbf{x}\|_1 \|\mathbf{y}\|_\infty. $
Furthermore, to equate the powers of the elements, we incorporated $ |z|^{-1} $ into the constant $ D $, keeping in mind that $ D $ is a generic constant and its value can be adjusted as needed.


Now, by considering the definition of $\tilde\bR_{k,k+1-l}$ in \eqref{eq: tildeR_def}, or specifically its non-zero elements, we define the $p^\tth$ delay power spectrum of node $(v,u)$ and its delay response at time $m$, respectively, as
\begin{equation}
    \tau^{vu}_{p, k}(\beta)=\frac 1{|T_k^{vu}|}\suml_{k^\prime\in T_k^{vu}} z_p^{|k+1-k^\prime|}(\beta)-1,
\end{equation}
\begin{equation}
    \kappa^{vu}_{m}(\beta)=\suml_p\suml_{k\mid T_k^{vu}\neq\emptyset}|z_p(\beta)|^{|k-m|}\tau^{vu}_{p, k}(\beta).
\end{equation}
Then, the delay response $\kappa$ is defined as the supremum of $\kappa^{vu}_{m}(\beta)$ over every time $m$ and link $(v,u)$ and $\beta$. Therefore, the inequality represented in \eqref{eq: 83} simplifies to
\begin{equation}
    \bar\lambda \leq D\kappa \lambda. 
\end{equation}

From \eqref{eq: 76_new}, we also have $\lambda \leq D\bar\lambda$. Combining these two inequalities leads us to 
\begin{equation} \label{eq: 93}
    \lambda \leq D\kappa\lambda. 
\end{equation}

Keep in mind that $D$ is a generic constant, and it relies on $C(\kappa)$. Specifically, when $C$ increases, $D$ decreases, but it will not go lower than $1/\kappa_0$.
If we have a delay response, $\kappa$, from one delay event, and it satisfies $\kappa < \kappa_0$, then there are possible values of $C$ where $D\kappa$ is less than 1. Based on \eqref{eq: 93}, this makes $\lambda=0$, which leads to $\tlbPsi_m = \bO$ and $\bLambda_{m-1}=\bO$. This gives us the full proof.

\newpage
\section*{Lemmas}
\begin{lemma}[Inverse-matrix decomposition] \label{lemma: matrix_decomposition}
let $\bF(z)$ be a square matrix of size $m$, and the set $\{ z_1, z_2, \ldots, z_P\}$ contains all the simple roots of $\det(\bF(z))=0$. Further, assume $z\bF_\beta^{-1}(z)$ has a finite limit as $z\to 0$.Then, the inverse matrix $\bF^{-1}(z)$ can be decomposed into 
\begin{equation}
    \bF^{-1}(z) = \suml_{p=1}^P \frac{\ba_p\bb_p^T}{z-z_p},
\end{equation} 
where
\[\ba_p \in \ker(\bF(z_p))\]
\[\bb_p \perp \text{image}(\bF(z_p)).\]

\begin{proof}
    From Cramer's rule for matrix inversion, considering that $\det\bF(z)=0$ has only simple roots, we decompose $\bF^{-1}(z)$ into:
    \begin{equation}
        \bF^{-1}(z) = \suml_{p=1}^P \frac{\bL_p}{z-z_p},
    \end{equation}
    where $\bL_p$ are some Known matrices. Since $\bF(z)$ is invertible for every element $\bb^\prime$ in the range space (image) of $\bF(z)$, for all $z$, the following condition should be hold
    \begin{equation}\label{eq:L_null}
     \bL_p\bb^\prime = \bzero,
    \end{equation}
    otherwise the inverse does not exist at $z=z_p$.

    From the fundamental theorem of linear algebra, specifically the Rank-Nullity theorem, the rank or the dimension of the range space of $\bF(z_p)$ is $m-1$ since the nullity, i.e. the dimension of the null space, of $\bF(z_p)$ is one. Therefore, the range space can be represented by only one vector orthogonal to this subspace, represented by $\bb_p$. 

    Since any vector in the range space of $\bF(z)$, containing $\bb^\prime$, is orthogonal to $\bb_p$, $\bL_p$ should be in the form of $\bL_p=\ba_p\bb_p^T$ to satisfy \eqref{eq:L_null}. By following the same procedure for $\bF(z)^T$, considering 
    $
    \det(\bF(z)^T) = \det(\bF(z)),
    $
    we observe that $\ba_p$ is an element in the null space of $\bF(z_p)$.  
\end{proof}
\end{lemma}

\begin{lemma}[Inversely-related roots] \label{lemma: inversly-related_roots}
    let matrix $\bF(z)$ be in the following form
    \begin{equation}
        \bF(z) = \left[  
        \begin{array}{cc}
            \bA & z^{-1}\bB^T \\
            z\bB & \bC
        \end{array}
        \right],
    \end{equation}
    with symmetric matrices $\bA,\bB$. Then, the roots of $\det(\bF(z))$ are \emph{inversely-related}, meaning that if $z_i$ is a root, then its inverse $z_i^{-1}$ is also a root. 

\begin{proof}
Given the configuration of $\bF(z)$, we have $\bF(z)^T = \bF(z^{-1})$. Because $\det\bF(z)^T = \det\bF(z^{-1}) = \det\bF(z)$, it follows that if $z_i$ serves as a root, its reciprocal $z_i^{-1}$ is also a root.
\end{proof}
\end{lemma}


\begin{lemma} \label{lem:integral}
[Contour integral over the unit circle]
Let  $f(z) = \frac{z^k}{z-z_p}$. The contour integral of $f(z)$ over the unit circle $\Gamma$ is given by
\begin{equation}\label{eq: solution_integral}
    \frac{1}{2\pi j}\oint_\Gamma \, \frac{z^k}{z-z_p} = 
    \begin{cases}
    z_p^k & |z_p| < 1, k\geq0 \\
    -z_p^k & |z_p| > 1, k<0 \\
    0 & \text{otherwise}.
    \end{cases}
\end{equation}

\begin{proof}
    To compute the integral, we consider the following four cases:
    \begin{itemize}
        \item $|z_p|>1, k\geq0$: The function is analytic inside the contour. Based on the Cauchy integral theorem, the integral equals $0.$
        \item $|z_p|<1, k\geq0$: $z=z_p$ is the only singular point of $f(z)$ inside the contour. Based on the residue theorem, the integral equals $z_p^k$. 
        \item $|z_p|>1, k<0$:
        $z=0$ is the only singularity of order $k$ of $f(z)$. The residue at $z=0$ is computed as
        \[
        \text{Res}[f(z), 0] = \lim_{z\to 0} \frac{1}{(k-1)!} \frac{d^{k-1}}{dz^{k-1}} \left[ \frac{1}{z-z_p} \right] = -z_p^{k}
        \]
        Based on the Residue Theorem, the solution of the integral is $ -z_p^k $.
        \item $|z_p|<1, k<0$: The function $ f(z) $ has two singularities at $ z=0 $ and $ z=z_p $. Based on the Residue Theorem, considering that the residues cancel each other out, the integral equals 0.
    \end{itemize}
Combining the above cases, we obtain the solution to the contour integral as given in \eqref{eq: solution_integral}.
\end{proof}
\end{lemma}

\end{document}